\documentclass{amsart}
\usepackage{amsmath, amsthm, amscd, amsfonts, amssymb, graphicx, xcolor}
\usepackage[bookmarksnumbered, colorlinks, plainpages]{hyperref}
\usepackage{cite,comment}
\usepackage{graphicx}
\usepackage{float}
\usepackage[none]{hyphenat}

\usepackage{enumerate}
\usepackage{xcolor, soul}

\newtheorem{theorem}{Theorem}[section]
\newtheorem{lemma}[theorem]{Lemma}

\newtheorem{corollary}[theorem]{Corollary}
\theoremstyle{definition}
\newtheorem{definition}[theorem]{Definition}

\theoremstyle{remark}
\newtheorem{remark}[theorem]{Remark}
\numberwithin{equation}{section}

\begin{document}
	\setcounter{page}{1}
	\centerline{}
	
	\centerline{}
	
	\title[Sharp inequalities for fractional spherical maximal operators on finite graphs]{Sharp inequalities for fractional spherical maximal operators on finite graphs}
	\author[S. Lu, Q. shen and F. Zhao]{Shengzhe Lu, Qizhen Shen and Fayou Zhao$^*$}
	\address[Shengzhe Lu]{Department of Mathematics, Shanghai University, Shanghai
		200444, P.R. China}
	\email{lushengzhe@shu.edu.cn}
	
	\address[Qizhen Shen]{Qian Weichang College, Shanghai University, Shanghai
		200444, P.R. China}
	\email{qizhen$_{-}$shen@163.com}
	
	\address[Fayou Zhao]{Newtouch Center for Mathematics of Shanghai University, Department of Mathematics,
		Shanghai University, Shanghai 200444, P.R. China}
	\email{fyzhao@shu.edu.cn}
	\thanks{The research was supported by National Natural Science Foundation of
		China (Grant No. 11971295).\\
		$\ast$ Corresponding author}
	
\begin{abstract}
This paper presents a systematic study of (quasi-)norm bounds for the fractional spherical maximal operator $M_{\alpha,G}^{\circ}$ on finite connected graphs $G$. We establish sharp estimates for the operator norm of $M_{\alpha,G}^{\circ}$ from $\ell^p(G)$ to $\ell^q(G)$ for $0 < p \leq q \leq \infty$, with explicit extremal values attained by complete graphs $K_n$ and star graphs $S_n$. We solve the corresponding extremal graph problem by showing that complete and star graphs are the unique maximizers of the operator norm in different parameter regimes. Additionally, we characterize graphs containing a universal vertex in terms of their operator norm. Finally, we analyze the spherical maximal operator $M^{\circ}_{G}$ on paths $P_n$, uncovering a clear divergence in the operator's behavior between odd and even lengths.
\end{abstract}

	\maketitle
	
\section{Introduction}

For any locally integrable function $f$ on $\mathbb{R}^n$, the fractional Hardy--Littlewood maximal operator $\mathcal{M}_{\alpha}$, with $0 \leq \alpha < n$, is defined by
\[
\mathcal{M}_{\alpha}f(x) = \sup_{r > 0} \frac{1}{|B(x,r)|^{1 - \alpha/n}} \int_{B(x,r)} |f(y)|  dy, \quad x \in \mathbb{R}^n.
\]
Similarly, for a suitable function $f$ on $\mathbb{R}^n$, the spherical maximal operator $\mathcal{M}^{\circ}_{\beta}$, with $0 \leq \beta < n-1$, is given by
\[
\mathcal{M}^{\circ}_{\beta}f(x) = \sup_{r > 0} \frac{1}{|S(x,r)|^{1 - \beta/(n-1)}} \int_{S(x,r)} |f(y)|  d\sigma(y), \quad x \in \mathbb{R}^n,
\]
where $B(x,r)$ denotes the open ball centered at $x$ with radius $r$, $S(x,r)$ is the sphere centered at $x$ with radius $r$, $d\sigma$ is the normalized rotation-invariant measure on the unit sphere in $\mathbb{R}^{n}$, and $|E|$ represents the $n$-dimensional Lebesgue measure of a measurable set $E \subset \mathbb{R}^n$. The special cases $\alpha = 0$ and $\beta = 0$ are denoted by $\mathcal{M} := \mathcal{M}_{0}$ and $\mathcal{M}^{\circ} := \mathcal{M}^{\circ}_{0}$, respectively.

The Hardy--Littlewood--Wiener theorem asserts that the Hardy--Littlewood maximal operator $\mathcal{M}$ is of weak type $(1, 1)$ and of strong type $(p, p)$ for all $1 < p \leq \infty$; see \cite{HL1930, WN1939}. For $0 \leq \alpha < n$ and $f \in L^1_{\mathrm{loc}}(\mathbb{R}^{n})$, the fractional integral operator of order $\alpha$, denoted $\mathcal{I}_{\alpha}$, is defined by
\[
\mathcal{I}_{\alpha}(f)(x) = \int_{\mathbb{R}^{n}} \frac{f(y)}{|x - y|^{n - \alpha}}  dy, \quad x \in \mathbb{R}^{n}.
\]
The strong type $(p, q)$ estimates for $\mathcal{I}_{\alpha}$ were first established by Hardy and Littlewood \cite{HL1928, HL1932} in the one-dimensional case, and later extended to higher dimensions by Sobolev \cite{S1938}, under the relation $1/p - 1/q = \alpha/n$ for $1 < p \leq q < \infty$. As for the endpoint estimate, the weak-type $\big(1, \frac{n}{n - \alpha}\big)$ bound was established by Zygmund \cite{Z1956} in the case $n = 1$, and later by Sobolev \cite{S1938} and Thorin \cite{T1948} for general dimensions. The classical technique for establishing the boundedness of the fractional integral operator $\mathcal{I}_{\alpha}$ was later provided by Hedberg \cite{Hed1972}. The boundedness of $\mathcal{M}_{\alpha}$ can be established by means of a pointwise inequality
\[
\mathcal{M}_{\alpha}f(x) \leq \mathcal{I}_{\alpha}(|f|)(x) \quad \text{for a.e. } x \in \mathbb{R}^n,
\]
as shown in \cite{CWW1993}; Muckenhoupt \cite{M1974} directly proved analogous results in the weighted setting.

The spherical maximal operator $\mathcal{M}^\circ$ is bounded on $L^p(\mathbb{R}^n)$ if and only if $ n/(n-1)<p\leq \infty$. This result was first obtained by Stein \cite{S1976} for $n \geq 3$, and by Bourgain \cite{B1986} for $n = 2$. At the endpoint $p = n/(n-1)$, Bourgain \cite{B1985} established a restricted weak-type inequality for $n \geq 3$, while the corresponding result fails for $n = 2$ as shown in \cite{STW2003}. For $n \geq 3$, $0 \leq \beta < n-1$, $n/(n-1) < p \leq \infty$, and $p \leq q < (n-1)p$ satisfying $1/p - 1/q = \beta/n$, Oberlin \cite[Theorem 2, p.146]{O1989} proved that $\mathcal{M}^{\circ}_{\beta}$ is of strong type $(p, q)$ using interpolation techniques; see \cite{CCG2002} for related results.

In this paper, we study norm estimates for the fractional spherical maximal operator  on finite connected graphs. This continues a classical line of research concerning maximal operators in Euclidean spaces, extending it to the discrete setting of graphs. We now introduce the necessary definitions and notation.
	
Let $G = (V_G, E_G)$ be a finite, undirected, connected graph with vertex set $V_G$ and edge set $E_G$. We endow the graph with the metric structure induced by the geodesic distance. Two vertices $u, v \in V_G$ are said to be adjacent, denoted $u \sim v$, if they are connected by an edge in $E_G$. For any $u, v \in V_G$, let $d_G(u, v)$ denote the number of edges in a shortest path connecting $u$ and $v$. Abusing notation, for a finite set $F$, we denote its cardinality by $|F|$. The degree of a vertex $u \in V_G$ is defined as $\deg(u) = |\{v \in V_G : d_G(u, v) = 1\}|$.
For convenience, we denote the vertex set by $V_G = \{1, 2, \ldots, n\}$ when $|V_{G}|=n$, let $\mathbb{R}_+ := [0, \infty)$, and let $\mathbb{N}_{+}:=\{1,2,\ldots\}$.

We also define three special families of graphs that will be of interest:
(i) The complete graph $K_n$: a graph with $n$ vertices such that every vertex has degree $n-1$;
(ii) The star graph $S_n$: a graph with $n$ vertices where one vertex (the center) has degree $n-1$, and all the others (the leaves) have degree $1$.
(iii) The path $P_n$: a graph with $n$ vertices such that two vertices have degree $1$ and all others have degree $2$.

The fractional spherical maximal operator on a graph $G$, denoted by $M_{\alpha,G}^{\circ}$ for $0 \leq \alpha < 1$, is defined for functions $f: V_G \to \mathbb{R}$ as
\[
M_{\alpha,G}^{\circ}f(x) = \sup_{r \geq 0} \frac{1}{|S(x,r)|^{1-\alpha}} \sum_{v \in S(x,r)} |f(v)|,
\]	
where $S(x,r) = \{ v \in V_G : d_G(v,x) = r \}$ denotes the sphere of radius $r$ centered at vertex $x$.
Correspondingly, the fractional Hardy--Littlewood maximal operator on the graph $G$, denoted by $M_{\alpha,G}$ (see \cite{KP1986}), is defined for functions $f: V_G \to \mathbb{R}$ by
\[
M_{\alpha,G}f(x) = \sup_{r \geq 0} \frac{1}{|B(x,r)|^{1-\alpha}} \sum_{v \in B(x,r)} |f(v)|,
\]	
where $B(x,r) =\{ v \in V_G : d_G(v,x) \leq r\}$ denotes the closed ball of radius $r$ centered at $x$. 
For notational simplicity, we adopt the conventions $M_{G}^{\circ} = M_{0,G}^{\circ}$ and $M_{G} = M_{0,G}$ throughout this work.
	
The study of averaging operators on finite connected graphs $G$, including the Hardy--Littlewood maximal operator, the fractional Hardy--Littlewood maximal operator, the spherical maximal operator, and the Hardy operator, has attracted considerable attention. A central objective in this line of research is the determination of sharp operator norm constants on graphs. Establishing such constants is of substantial interest, as improvements in this direction can refine quantitative estimates in various applications of the Hardy--Littlewood inequality, including quantitative versions of Rademacher's theorem for Lipschitz functions \cite{B1985-2}. The investigation is also motivated by discretization results for the Hardy--Littlewood maximal function in $\mathbb{R}$, particularly those concerning measures given by positive linear combinations of Dirac delta functions \cite{Mel2003,MS1992}, which are closely related to the setting of linear trees, see also \cite{ST2016}.

Research has demonstrated that graphs possessing richer geometric structures often admit improved estimates for averaging operators, and that these estimates can, in certain extremal cases, serve to characterize the graph itself. Initial insights into the connections between graph geometry and the associated maximal operator $M_G$ were provided by Soria and Tradacete \cite{ST2016}. They established sharp bounds for the $\ell^p$-(quasi)norm of $M_G$ over general graphs $G$ for $0 < p \leq 1$, identifying the complete graph $K_n$ and the star graph $S_n$ as the extremal graphs attaining the lower and upper bounds, respectively. The proof of sharp norms in this regime relies heavily on the concavity of the function $x \to x^p$, with the (quasi-)norm $\Vert M_G\Vert_{\ell^p\to \ell^p}$ being attained by a suitable Kronecker delta function. In contrast, the case $p > 1$ presents substantial challenges, as noted in~\cite[Remark 2.8]{ST2016}.
	
The first advances for $p > 1$ were made by Gonz\'alez-Riquelme and Madrid in \cite{CJ2021}, where the exact values of $\|M_{K_{n}}\|_{\ell^2\to \ell^2}$ and $\|M_{S_{n}}\|_{\ell^2\to \ell^2}$ were obtained, along with their corresponding extremizers. Subsequent work in \cite{CJ2022} provided a complete characterization of the extremizers for $\|M_{S_{n}}\|_{\ell^p\to \ell^p}$ for all $p \in (1,2]$, together with partial characterizations for $p > 2$ and an analogous result for the extremizers of $\|M_{K_{n}}\|_{\ell^p\to \ell^p}$ when $p > 1$.
	
Further extending this line of inquiry, Xiang and Zhao \cite{XZ2025} derived both lower and upper bounds for the optimal constant $C_{G,p,q}$ in the inequality $\Vert M_{\alpha,G} f \Vert_q \leq C_{G,p,q} \Vert f \Vert_p$ in the range $0 < p \leq q \leq 1$, again identifying $K_n$ and $S_n$ as extremal graphs attaining the two bounds. For $1 < p < \infty$, they also proved the existence of extremizers for the norms $\Vert M_{\alpha, K_n} \Vert_{\ell^p \to \ell^p}$ and $\Vert M_{\alpha, S_n} \Vert_{\ell^p \to \ell^p}$. Recently, Wang and Zhao \cite{WZ2025} established analogous optimal results for the $\ell^p$-(quasi)norm of the Hardy operator on finite graphs, including the cases where $G$ is a complete graph $K_n$, a star graph $S_n$, or a path $P_n$.
	
In a related direction, Zhang and Liu \cite{LZ20212} investigated the spherical maximal operator $M_{G}^{\circ}$ on finite graphs, establishing strong-type and weak-type $(p,p)$ estimates for all $0<p<\infty$, along with restricted type estimates on complete and star graphs. Hussain et al. \cite{HZFS2024} further characterized the extremal graphs realizing these bounds for $0<p\leq 1$.
	
Motivated by these foundational advances in the theory of averaging operators on graphs, this paper undertakes a systematic investigation of (quasi-)norm estimates for their fractional variants in the finite graph setting. It is readily observed that for any connected finite graph $G$ with $|V_G|=n$ and $0 \leq \alpha<1$, the fractional spherical maximal operator $M^{\circ}_{\alpha,G}$ and the fractional Hardy--Littlewood maximal operator $M_{\alpha,G}$ satisfy the following pointwise relationship:
\[
n^{-\alpha}M_{\alpha,G}f(x) \leq M^{\circ}_{\alpha,G}f(x) \leq n^{1-\alpha}M_{\alpha,G}f(x).
\]
While this inequality directly implies the boundedness of  $M^{\circ}_{\alpha,G}$  via known results for $M_{\alpha,G}$ in \cite{XZ2025}, obtaining sharp norm bounds for the spherical operator $M^{\circ}_{\alpha,G}$ on various function spaces and graph structures remains a significant and non-trivial challenge.

The remainder of this article is structured as follows. In Section \ref{s2}, we determine the (quasi-)norm $\left \Vert M_{\alpha,G}^{\circ} \right \Vert_{\ell^p \to \ell^q}$ for $0 < p \leq q \leq \infty$, with particular emphasis on the extremal cases where $G$ is the complete graph $K_n$ or the star graph $S_n$. In Section \ref{s3}, we shift our focus to the problem of finding the largest possible operator norm $\left\| M_{\alpha,G}^{\circ} \right\|_{\ell^p\to \ell^q}$ among all connected graphs with $n$ vertices. We completely resolve this extremal problem for a wide range of parameters. Our main results prove that the complete graph $K_n$ and the star graph $S_n$ are the unique extremal graphs maximizing the norm in two distinct parameter regimes. Furthermore, for $0 \leq \alpha < 1$ and $1/\alpha \leq p \leq \infty$, we provide a characterization of graphs containing a universal vertex based on the value of the operator norm $\| M^{\circ}_{\alpha,G} \|_{\ell^p \to \ell^\infty}$. These extremal graph results highlight the critical role of graph diameter and vertex centrality in determining the behaviour of $M_{\alpha,G}^{\circ}$.
In the last Section, we turn to the analysis of the spherical maximal operator $M^{\circ}_{G}$ on the path $P_n$. The behavior of this operator differs significantly depending on the parity of the path, presenting substantial technical challenges.
	
\section{Optimal estimates for $\left \Vert M_{\alpha,G}^{\circ} \right \Vert_{\ell^p\to \ell^q}$ } \label{s2}
This section is devoted to establishing the precise values and sharp estimates of the operator norm $\left\| M_{\alpha,G}^{\circ} \right\|_{\ell^p \to \ell^q}$ for finite connected graphs. We consider three fundamental graph structures: general connected graphs, complete graphs $K_n$, and star graphs $S_n$. For each graph type and for various ranges of the parameters $0 \leq \alpha < 1$ and $0 < p\leq q \leq \infty$, we provide either exact computations or tight upper and lower bounds of the norm.

\begin{definition}
Let $G=(V_G, E_G)$ be a graph, define
\begin{align*}
\left\|f\right\|_p:=\left\|f\right\|_{\ell^{p}(G)}=\Big(\sum\limits_{v\in V_G}|f(v)|^p\Big)^{1/p} \ \text{and}\
\left\|f\right\|_{\ell^{p,\infty}(G)}:=\underset{\lambda>0}{{\rm sup}}\lambda|\{v\in V_G:f(v)>\lambda\}|^{1/p}			
\end{align*} for $0<p<\infty$,
and $\left\|f\right\|_{\infty}:=\left\|f\right\|_{\ell^\infty}=\underset{v\in V_G}{\sup}|f(v)|$. The $\ell^p \to \ell^q$-norm (quasi-norm in the range $0<q<1$) of an operator $T$, acting on functions defined over a  graph $G$, is given by
\[
\left\|T\right\|_{\ell^p\to \ell^q}:=\sup_{\left\|f\right\|_{\ell^p(G)}\neq 0}\frac{\left\|Tf\right\|_{\ell^q(G)}}{\left\|f\right\|_{\ell^p(G)}}, \quad 0<p,q\leq \infty.
\]
And the $\ell^p\to \ell^{q,\infty}$-norm (quasi-norm in the range $0<q<1$) of an operator $T$ is given by
\[
\left\|T\right\|_{\ell^p\to \ell^{q,\infty}}:=\sup_{\left\|f\right\|_{\ell^p(G)}\ne 0}\frac{\left\|Tf\right\|_{\ell^{q,\infty}(G)}}{\left\|f\right\|_{\ell^p(G)}}, \quad 0<p,q\leq \infty.
\]		
\end{definition}
Before we start our work, we first give an important lemma \cite{XZ2025}.
\begin{lemma}\label{lem2-2}
Let $G=(V_G, E_G)$ be a finite connected graph with $|V_G|=n\geq 2$ vertices, and let $T:\ell^{p}(G)\to\ell^{q}(G)$ be a sublinear operator  with $0<p\leq q\leq1$. Then
$$\Vert T\Vert_{\ell^{p}\to\ell^{q}}=\max_{k\in V_{G}}\Vert T\delta_k\Vert_q,$$
where $\delta_k$ is the Kronecker delta function defined as
\begin{equation*}
\delta_k(i)=
\begin{cases}
	1 & {i=k,}\\
	0 & {i\neq k.}
\end{cases}	
\end{equation*}
\end{lemma}

We now establish the foundational estimates for the operator (quasi-)norm of the fractional spherical maximal operator on a general finite connected graph. The following theorem provides both upper and lower bounds, which are applicable across various ranges of the exponents $p$ and $q$.
\begin{theorem}\label{thm2-3}
Let $G=(V_G, E_G)$ be a finite connected graph with $|V_G|=n\geq 3$ and $0\leq \alpha<1$.

{\rm(i)} If $0<p\leq q\leq1$, then
\[
\left(1+(n-1)^{q(\alpha-1)+1}\right)^{1/q}\leq \left \| M_{\alpha,G}^{\circ} \right \|_{\ell^p\to \ell^q}\leq n^{1/q}.
\]

{\rm(ii)} If $1<p\leq q<\infty$, then
\[
\left(1+(n-1)^{q(\alpha-1)+1}\right)^{1/q}\leq \left \| M_{\alpha,G}^{\circ} \right \|_{\ell^p\to \ell^q}\leq \left\{\begin{matrix}
	n^{1/q} & p\leq1/\alpha,\\[0.2cm]
	n^{1/q}(n-1)^{\alpha-{1/p}}&p>1/\alpha
\end{matrix}\right.
\]
for $0<\alpha<1$, and
\[
\left(1+(n-1)^{1-q}\right)^{1/q}\leq \left \| M_{G}^{\circ} \right \|_{\ell^p\to \ell^q}\leq n^{1/q}
\]
for $\alpha=0$.

{\rm(iii)} If $0<p\leq1< q<\infty$, then
\[
\left(1+(n-1)^{q(\alpha-1)+1}\right)^{1/q}\leq \left \| M_{\alpha,G}^{\circ} \right \|_{\ell^p\to \ell^q}\leq n^{1/q}.
\]

{\rm(iv)} If $q=\infty$, then
\[
\left \| M_{\alpha,G}^{\circ} \right \|_{\ell^p\to \ell^\infty}= 1
\]
for $0<\alpha<1$, $0<p\leq {1}/{\alpha}$, and
\[
1\leq \left \| M_{\alpha,G}^{\circ} \right \|_{\ell^p\to \ell^\infty}\leq (n-1)^{\alpha-{1/p}}
\]
for $0<\alpha<1$, ${1}/{\alpha}<p< \infty$, and
\[
\left \| M_{G}^{\circ} \right \|_{\ell^p\to \ell^\infty}= 1
\]
for $\alpha=0$, $0<p<\infty$.

{\rm(v)} If $p=q=\infty$, then
\[1\leq \left \| M_{\alpha,G}^{\circ} \right \|_{\ell^\infty\to \ell^\infty}\leq (n-1)^{\alpha}.\]

{\rm(vi)} If $G$ is a finite connected graph with $|V_G|=2$, then
\begin{align*}
	\left \| M_{\alpha,G}^{\circ} \right \|_{\ell^p\to \ell^q}= \left\{\begin{matrix}
		2^{1/q} & 0<p\leq q<\infty, \\
		1 & 0<p\leq \infty,\  q=\infty.
	\end{matrix}\right.
\end{align*}	
\end{theorem}

\begin{proof}	
Noting that for any $v\in V_G$, we have
\begin{align}\label{equ2-1}
M_{\alpha,G}^{\circ}f(v)&=\underset{r=0,\ldots,n-1}{\max}\frac{1}{|S(v,r)|^{1-\alpha}}\sum_{x\in S(v,r)}|f(x)|.
\end{align}
Throughout this paper, if there exists some $r_0 \in \{0, \ldots, n-1\}$ such that $|S(v, r_0)| = 0$, then for any $0 \leq \alpha < 1$, we define
\[
\frac{1}{|S(v, r_0)|^{1-\alpha}} \sum_{x \in S(v, r_0)} |f(x)| = 0.
\]
Accordingly, we henceforth assume $S(v, r)$ is non-empty and thus $1 \leq |S(v, r)| \leq n - 1$.

For any $0<p\leq \infty$, fix $i\in V_G$ and take $f=\delta_i$, then $\Vert\delta_i\Vert_{p}=1$. And for any $0<q<\infty$, we have
\begin{align*}
	\left \| M_{\alpha,G}^{\circ}\delta_i\right \|_{q}&=\bigg( (M_{\alpha,G}^{\circ}\delta_i(i))^{q}+\sum_{v \in V_G\setminus \{i\}}( M_{\alpha,G}^{\circ}\delta_i(v))^q\bigg)^{1/q}\\
	&=\bigg(1+\sum_{v\in V_G\setminus \{i\}}\bigg(\max_{r=1,\ldots,n-1}\frac{1}{|S(v,r)|^{1-\alpha}}\sum_{x\in S(v,r)}\delta_i(x)\bigg)^q\bigg)^{1/q}.
\end{align*}

Since $1\leq|S(v,r)|\leq n-1$ for all $ v\in V_G\setminus \{i\}$, it holds for any $0\leq \alpha< 1$ that
\[
\left(1+(n-1)^{q(\alpha-1)+1}\right)^{1/q}\leq \left \| M_{\alpha,G}^{\circ}\delta_{i} \right \|_{q}\leq n^{1/q}.
\]

{\rm(i)} The desired result follows immediately from Lemma \ref{lem2-2} for the range $0 < p \leq q \leq 1$.

{\rm(ii)}  It follows from \eqref{equ2-1} and the finiteness of the set of radii that for any vertex $v\in V_{G}$,  there exists an $r_v\in \{0,1,\ldots,n-1\}$ such that
\begin{align}\label{equ2-2}
M_{\alpha,G}^{\circ}f(v)&=|S(v,r_v)|^{\alpha-1}\sum_{x\in S(v,r_v)} |f(x)|.	
\end{align}
For $1<p<\infty$ and $0<\alpha<1$, by H\"{o}lder's inequality, we have
\begin{align}\label{p>1}
\notag M_{\alpha,G}^{\circ}f(v)&\leq |S(v,r_v)|^{\alpha-{1/p}}\bigg(\sum_{x\in S(v,r_v)} |f(x)|^p\bigg)^{1/p} \\
&\leq \left\{\begin{matrix}
	\left\|f\right\|_p & 1<p\leq1/\alpha,\\[0.2cm]
	(n-1)^{\alpha-{1/p}}\left\|f\right\|_p& p>1/\alpha .
\end{matrix}\right.	
\end{align}
This leads to
\begin{align*}
\left \| M_{\alpha,G}^{\circ}\right \|_{\ell^p\to \ell^q}\leq \left\{\begin{matrix}
	n^{1/q} & 1<p\leq1/\alpha,\\[0.2cm]
	(n-1)^{\alpha-{1/p}}n^{1/q}&p>1/\alpha .
\end{matrix}\right.	
\end{align*}

For $1<p<\infty$ and $\alpha=0$, by H\"{o}lder's inequality, we have
\begin{align}\label{equ77-2}
	M_{G}^{\circ}f(v)\leq |S(v,r_v)|^{-{1/p}}\bigg(\sum_{x\in S(v,r_v)} |f(x)|^p\bigg)^{1/p}\leq \left\|f\right\|_{p},
\end{align}
which leads to
\[
\left \| M_{G}^{\circ}\right \|_{\ell^p\to \ell^q}\leq n^{1/q}.
\]

{\rm(iii)} For $0<p\leq 1$, it follows from \eqref{equ2-2} that for $0\leq \alpha<1$, we have
\begin{align}\label{p<1}
M_{\alpha,G}^{\circ}f(v)\leq  |S(v,r_v)|^{\alpha-1}\bigg(\sum_{x\in S(v,r_v)} |f(x)|^p\bigg)^{1/p}  \leq \left\|f\right\|_p.	
\end{align}
Then we have
\[
\left\|M_{\alpha,G}^{\circ}\right\|_{\ell^p\to \ell^q}\leq n^{1/q}.
\]

{\rm(iv)} It is easy to see that
\begin{align*}
\left\|M_{\alpha,G}^{\circ}f\right\|_\infty&=\max_{v\in V_G}\bigg\{\max_{r=0,\ldots,n-1}\frac{1}{|S(v,r)|^{1-\alpha}}\sum_{x\in S(v,r)}|f(x)| \bigg\}.	
\end{align*}
It follows from (ii) and (iii) that for $0<\alpha<1$, and any vertex $v\in V_G$, there exists an $r_v\in \{0,1,\ldots,n-1\}$ such that 
\begin{align*}
M_{\alpha,G}^{\circ}f(v)\leq \left\{\begin{matrix}
	|S(v,r_v)|^{\alpha-1} \left\|f\right\|_p & 0<p\leq1,\\[0.2cm]
	|S(v,r_v)|^{\alpha-{1/p}} \left\|f\right\|_p & 1<p<\infty.
\end{matrix}\right.	
\end{align*}
Thus we have
\begin{align*}
\left\|M_{\alpha,G}^{\circ}f\right\|_\infty\leq \left\{\begin{matrix}
	\left\|f\right\|_p  & 0<p\leq {1}/{\alpha},\\[0.2cm]
	(n-1)^{\alpha-{1/p}}  \left\|f\right\|_p & {1}/{\alpha}<p<\infty.
\end{matrix}\right.	
\end{align*}

It is straightforward to verify that for $\alpha=0$, the inequality
\[
\left\|M_{G}^{\circ}f\right\|_\infty\leq \left\|f\right\|_p
\]
holds.

The lower bound follows from taking $f = \delta_v$ for any $v \in V_G$ when $0 \leq \alpha < 1$.

{\rm(v)} The upper bound can be obtained by \eqref{equ2-1}, while the lower bound is found by setting $f = \delta_v$ for any $v \in V_G$.

{\rm(vi)} Without loss of generality, we assume $V_G=\{a,b\}$, $E_G=\{a\sim b\}$. For any $0<p\leq \infty$, we have
\[
M_{\alpha,G}^{\circ}f(a)=M_{\alpha,G}^{\circ}f(b)=\max\{|f(a)|,|f(b)|\}\leq \left\|f\right\|_p.
\]
Then for $0<q<\infty$, we have
\begin{align*}
\left\|M_{\alpha,G}^{\circ}\right\|_{\ell^p\to \ell^q}=\sup_{\left\|f\right\|_{\ell^p(G)}\neq 0} \frac{\left((M_{\alpha,G}^{\circ}f(a))^q+(M_{\alpha,G}^{\circ}f(b))^q\right)^{1/q}}{\left\|f\right\|_p}
\leq 2^{1/q}.	
\end{align*}
For $q=\infty$, we have
\begin{align*}
\left\|M_{\alpha,G}^{\circ}\right\|_{\ell^p\to \ell^\infty}= \sup_{\left\|f\right\|_{\ell^p(G)}\neq 0} \frac{\left\|M_{\alpha,G}^{\circ}f\right\|_\infty}{\left\|f\right\|_p} \leq 1.	
\end{align*}

To show the constants are sharp, we take $f=\delta_a$, then
\begin{align*}
\left\|M_{\alpha,G}^{\circ}\delta_a\right\|_q=\left\{\begin{matrix}
	2^{1/q} & 0<q<\infty,\\
	1 & q=\infty.
\end{matrix}\right.	
\end{align*}
This concludes the proof.
\end{proof}		
	
We now turn to the case of the complete graph $K_n$. Due to its high symmetry and uniform structure, we are able to derive significantly sharper bounds for the operator norm of
$M_{\alpha,K_n}^{\circ}$, and in many cases, establish its exact value. The following theorem provides a complete characterization.

\begin{theorem}\label{thm2-4} Let $0\leq \alpha<1$ and $n\geq 3$.
	
{\rm(i)} If $0<p\leq q\leq 1$, then
\[
\left \| M_{\alpha,K_n}^{\circ} \right \|_{\ell^p\to \ell^q}=\left(1+(n-1)^{q(\alpha-1)+1}\right)^{1/q}.
\]

{\rm(ii)}  If $1<p\leq q<\infty$, then
\begin{align*}
	\left(1+(n-1)^{q(\alpha-1)+1}\right)^{1/q}\leq \left \| M_{\alpha,K_n}^{\circ} \right \|_{\ell^p\to \ell^q}\leq \left\{\begin{matrix}
		\left(1+(n-1)^{q(\alpha-{1/p})+1}\right)^{1/q} & 1<p\leq1/\alpha,\\[0.2cm]
		n^{1/q}(n-1)^{\alpha-{1}/{p}} &p>1/\alpha
	\end{matrix}\right.
\end{align*}
for $0<\alpha<1$, and
\begin{align*}
	\left(1+(n-1)^{1-q}\right)^{1/q}\leq \left \| M_{K_n}^{\circ} \right \|_{\ell^p\to \ell^q}\leq \left(1+(n-1)^{1-q/p}\right)^{1/q}
\end{align*}
for $\alpha=0$.

{\rm(iii)}  If $0<p\leq1< q<\infty$, then
\[
\left \| M_{\alpha,K_n}^{\circ} \right \|_{\ell^p\to \ell^q}=\left(1+(n-1)^{q(\alpha-1)+1}\right)^{1/q}.
\]

{\rm(iv)} If $0<p<\infty$ and $q=\infty$, then
\begin{align*}
	\left \| M_{\alpha,K_n}^{\circ} \right \|_{\ell^p\to \ell^\infty}
	=\left\{\begin{matrix}
		1 & 0<p\leq {1}/{\alpha}, \\[0.2cm]
		(n-1)^{\alpha-{1/p}}  &1/{\alpha}<p<\infty
	\end{matrix}\right.
\end{align*}
for $0<\alpha<1$, and
\begin{align*}
	\left \| M_{K_n}^{\circ} \right \|_{\ell^p\to \ell^\infty}=1
\end{align*}
for $\alpha=0$.

{\rm(v)} If $p=q=\infty$, then
\[
\left \| M_{\alpha,K_n}^{\circ} \right \|_{\ell^\infty\to \ell^\infty}=(n-1)^\alpha.
\]	
\end{theorem}

\begin{proof}
If $G=K_n$, then for any $v\in V_{K_n}$, we have
\[
M_{\alpha,K_n}^{\circ} f(v)=\max_{r=0,1}\frac{1}{|S(v,r)|^{1-\alpha}}\sum_{x\in S(v,r)}|f(x)|.
\]
Fix $i \in V_{K_n}$ and set $f = \delta_i$. Note that $\Vert\delta_i\Vert_p = 1$ for any $0 < p \leq \infty$. Then, for any $0 < q < \infty$,
\begin{align*}
	\left \| M_{\alpha,K_n}^{\circ}\delta_i\right \|_{q}&=\bigg( (M_{\alpha,K_n}^{\circ}\delta_i(i))^{q}+\sum_{v \in V_{K_n}\setminus \{i\}}( M_{\alpha,K_n}^{\circ}\delta_i(v))^q\bigg)^{1/q}\\
	&=\bigg(1+\sum_{v\in V_{K_n}\setminus \{i\}}\bigg(\max_{r=0,1}\frac{1}{|S(v,r)|^{1-\alpha}}\sum_{x\in S(v,r)}\delta_i(x)\bigg)^q\bigg)^{1/q}.
\end{align*}

By a direct calculation, we have that for $v\in V_{K_n}\setminus\{ i\}$,
\[
M_{\alpha,K_n}^{\circ}\delta_i(v)=\max_{r=0,1}\frac{1}{|S(v,r)|^{1-\alpha}}\sum_{x\in S(v,r)}\delta_i(x)=\frac{1}{|S(v,1)|^{1-\alpha}}\sum_{x\in S(v,1)}\delta_i(x)={(n-1)}^{\alpha-1},
\]
and for $v=i$,
\[
M_{\alpha,K_n}^{\circ}\delta_i(i)=\max_{r=0,1}\frac{1}{|S(i,r)|^{1-\alpha}}\sum_{x\in S(i,r)}\delta_i(x)=1.
\]
Thus we have for any $0<q<\infty$,
\begin{align}\label{equ2-3}
	\left \| M_{\alpha,K_n}^{\circ}\delta_i\right \|_{q}=\left(1+(n-1)^{q(\alpha-1)+1}\right)^{1/q}.
\end{align}

(i) By Lemma $\ref{lem2-2}$, for $0<p\leq q \leq 1$, we have $$\left \| M_{\alpha,K_n}^{\circ} \right \|_{\ell^p\to \ell^q}=\left(1+(n-1)^{q(\alpha-1)+1}\right)^{1/q},$$
which proves ${\rm(i)}$.

(ii) For $0<\alpha<1$, $1<p\leq q <\infty$, it is straight forward to check that
\[
\left \| M_{\alpha,K_n}^{\circ}f\right \|_{q}=\bigg(\sum_{i=1}^{n}\max\bigg\{|f(i)|^q, \bigg((n-1)^{\alpha-1}\sum_{v\in V_{K_n}\setminus\{i\}}|f(v)|\bigg)^q\bigg\}\bigg)^{1/q}.
\]
Applying H\"{o}lder's inequality, we have
\[\sum_{v\in V_{K_n}\setminus\{i\}}|f(v)|\leq  (n-1)^{1-{1/p}}\Vert f\Vert_{p}.\]
Therefore
$$\left\|M_{\alpha,K_n}^{\circ}f\right\|_q\leq\left(\sum_{i=1}^{n}\max\left\{|f(i)|^q,(n-1)^{q(\alpha-{1}/{p})}\Vert f\Vert_{p}^{q}\right\}\right)^{{1}/{q}}.$$

Assume $\Vert f\Vert_{p}=1$. We consider the following two cases.

{\bf{Case 1.}} If  $|f(i)|\leq (n-1)^{\alpha-{1}/{p}}$ for all $1\leq i\leq n$, then
\[
\left\|M_{\alpha,K_n}^{\circ}\right\|_{\ell^p\to \ell^q}=\sup_{\left\|f\right\|_p=1}\left\|M_{\alpha,K_n}^{\circ}f\right\|_q\leq n^{1/q}(n-1)^{\alpha-{1}/{p}}.
\]

{\bf{Case 2.}} If there exist one or more indices $i \in \{1,2,\ldots,n\}$ for which $ |f(i)| > (n-1)^{\alpha-{1}/{p}}$, then for $1<p\leq q<\infty$, we have 
\begin{align*}
	\left\|M_{\alpha,K_n}^{\circ}\right\|_{\ell^p\to \ell^q}
	&\leq \sup_{\left\|f\right\|_p=1} \Bigg(\sum_{\{i: |f(i)|> (n-1)^{\alpha-{1}/{p}}\}}|f(i)|^{q}+\sum_{\{j: |f(j)|\leq (n-1)^{\alpha-{1}/{p}}\}}(n-1)^{q(\alpha-{1/p})}\Bigg)^{{1}/{q}}\\
	&\leq \left(1+(n-1)^{q(\alpha-{1/p})+1}\right)^{{1}/{q}}.
\end{align*}

Combining Cases 1 and 2 leads to
\begin{align*}\left\|M_{\alpha,K_n}^{\circ}\right\|_{\ell^p\to \ell^q}&\leq \max\left\{n^{1/q}(n-1)^{\alpha-{1}/{p}}, \left(1+(n-1)^{q(\alpha-{1/p})+1}\right)^{1/q}\right\}\\
	&= \left\{\begin{matrix}
		(1+(n-1)^{q(\alpha-{1/p})+1})^{1/q} &p\leq 1/\alpha, \\[0.2cm]
		n^{1/q}(n-1)^{\alpha-{1}/{p}} &p>1/\alpha.
	\end{matrix}\right.
\end{align*}

The verification for the case $\alpha=0$ follows a similar method, and we omit the details here.

(iii) For $0<p\leq 1$, we have
\[
\left\|M_{\alpha,K_n}^{\circ}f\right\|_q\leq\left(\sum_{i=1}^{n}\max\left\{|f(i)|^q,(n-1)^{q(\alpha-1)}\left\|f\right\|_p^q\right\}\right)^{{1}/{q}}.
\]
Assume $\Vert f\Vert_{p}=1$. Following the similar estimates as that of (ii) above, we have
\begin{align*}
	\left\|M_{\alpha,K_n}^{\circ}\right\|_{\ell^p\to \ell^q}&\leq \max\left\{n^{1/q}(n-1)^{\alpha-1},\left(1+(n-1)^{q(\alpha-1)+1}\right)^{{1}/{q}}\right\}\\
	&=\left(1+(n-1)^{q(\alpha-1)+1}\right)^{{1}/{q}}.
\end{align*}

Then we obtain that
\[
\left\|M_{\alpha,K_n}^{\circ}\right\|_{\ell^p\to \ell^q}=\left(1+(n-1)^{q(\alpha-1)+1}\right)^{{1}/{q}}.
\]
with the help of \eqref{equ2-3}.

(iv) We only need to deal with the case $1/\alpha<p<\infty$ for $0<\alpha<1$. A direct computation gives that
\begin{align*}
	\left\|M_{\alpha,K_n}^{\circ}f\right\|_{\infty}&=\max_{v\in V_{K_n}}\max\bigg\{|f(v)|,\frac{1}{(n-1)^{1-\alpha}}\sum_{x\in V_{K_n}\setminus\{v\}}|f(x)|\bigg\}.
\end{align*}

Note that $|f(v)|\leq \left\|f\right\|_{p}$ for any $v\in V_{K_n}$. By H\"{o}lder's inequality and $1/\alpha<p<\infty$, we have
\begin{align}\label{equ2-4}
	\Vert M_{\alpha,K_n}^{\circ}\Vert_{\ell^p\to \ell^{\infty}} \leq  (n-1)^{\alpha-{1/p}}.
\end{align}

Taking $f=1-\delta_{i}$ for any $i\in V_{K_n}$, we have
$\left\|f\right\|_p=(n-1)^{{1}/{p}}$ and
\begin{align*}
	\left\|M_{\alpha,K_n}^{\circ}(1-\delta_i)\right\|_{\infty}
    =\max_{v\in V_{K_n}}\max\Big\{|(1-\delta_i)(v)|,\frac{1}{(n-1)^{1-\alpha}}\sum_{x\in V_{K_n}\setminus\{v\}}|(1-\delta_i)(x)|\Big\}
	=(n-1)^{\alpha}.
\end{align*}
Then
\[
\left\|M_{\alpha,K_n}^{\circ}\right\|_{\ell^p\to\ell^\infty}\geq \frac{\left\|M_{\alpha,K_n}^{\circ}f\right\|_{\infty}}{\left\|f\right\|_{p}}=(n-1)^{\alpha-{1/p}}.
\]

Thus for $1/\alpha<p<\infty$, we have
\begin{align}\label{equ2-5}
	\left\|M_{\alpha,K_n}^{\circ}\right\|_{\ell^p\to\ell^\infty}\geq (n-1)^{\alpha-{1/p}}.
\end{align}

Combining $\eqref{equ2-4}$ and $\eqref{equ2-5}$, we have
\[\left \| M_{\alpha,K_n}^{\circ} \right \|_{\ell^p\to \ell^\infty}=(n-1)^{\alpha-{1}/{p}}\]
for ${1}/{\alpha}<p<\infty$.

(v) For $p=q=\infty$, we take $f=1-\delta_i$ for any fixed vertex $i\in V_{K_n}$.
It is easy to see that $\left\|1-\delta_i\right\|_{\infty}=1$,
$M_{\alpha,K_n}^{\circ}(1-\delta_i)(i)=(n-1)^{\alpha}$ and
\begin{align*}
	M_{\alpha,K_n}^{\circ}(1-\delta_i)(v)&=\max\bigg\{1,\frac{n-2}{(n-1)^{1-\alpha}}\bigg\}
\end{align*} for $v\in V_G \setminus\{i\}$. This implies that
\[
\left\|M_{\alpha,K_n}^{\circ}(1-\delta_i)\right\|_{\infty}=(n-1)^\alpha.
\]

We now turn to estimating the upper bound. Observe that
\begin{align*}
	\left\|M_{\alpha,K_n}^{\circ}f\right\|_{\infty}
	&=\max_{v\in V_{K_n}}\max\Big\{|f(v)|,\frac{1}{(n-1)^{1-\alpha}}\sum_{x\in V_{K_n}\setminus\{v\}}|f(x)|\Big\}\\
	&\leq \max_{v\in V_{K_n}}\max\left\{\left\|f\right\|_{\infty},(n-1)^{\alpha}\left\|f\right\|_{\infty}\right\}.
\end{align*}

Hence we have
\[
\left\|M_{\alpha,K_n}^{\circ}\right\|_{\ell^\infty\to\ell^\infty}= (n-1)^\alpha.
\]	
\end{proof}

We now turn to the star graph $S_n$, which exhibits a distinct asymmetric structure in contrast to the complete graph.
In the following theorem, we establish estimates for the operator norm of $M_{\alpha,S_n}^{\circ} $.
The differing roles of the central vertex and the leaves play a crucial part in the analysis.

\begin{theorem}\label{thm2-6}
Let $0\leq \alpha<1$,  $n\geq 3$ and $0<p\leq q\leq \infty$.

{\rm(i)} If $0<p\leq q\leq 1$, then
$$\left \| M_{\alpha,S_n}^{\circ} \right \|_{\ell^p\to \ell^q}=n^{1/q}.$$

{\rm(ii)} If $1<p\leq q<\infty$, then
\begin{align*}
	\left \| M_{\alpha,S_n}^{\circ} \right \|_{\ell^p\to \ell^q}=n^{1/q}
\end{align*}
for $0<\alpha<1$, $1<p\leq1/\alpha$, and
\begin{align*}
	n^{1/q}\leq \left \| M_{\alpha,S_n}^{\circ} \right \|_{\ell^p\to \ell^q}\leq \left((n-1)^{q(\alpha-{1/p})}+(n-1)(n-2)^{q(\alpha-{1/p})}\right)^{1/q}
\end{align*}
for $0<\alpha<1$, $ 1/\alpha<p<\infty$, and
\begin{align*}
	\left \| M_{S_n}^{\circ} \right \|_{\ell^p\to \ell^q}=n^{1/q}
\end{align*}
for $\alpha=0$.

{\rm(iii)} If $0<p\leq 1<q<\infty$, then
$$\left \| M_{\alpha,S_n}^{\circ} \right \|_{\ell^p\to \ell^q}=n^{1/q}.$$

{\rm(iv)} If $0<p<\infty$ and $q=\infty$, then
\begin{align*}
	\left \| M_{\alpha,S_n}^{\circ} \right \|_{\ell^p\to \ell^\infty}
	=\left\{\begin{matrix}
		1 & 0<p\leq {1}/{\alpha},\\
		(n-1)^{\alpha-{1/p}}  &1/{\alpha}<p<\infty
	\end{matrix}\right.
\end{align*}
for $0<\alpha<1$, and
\begin{align*}
	\left \| M_{S_n}^{\circ} \right \|_{\ell^p\to \ell^\infty}=1
\end{align*}
for $\alpha=0$.

{\rm(v)} If $p=q=\infty$, then
$$\left \| M_{\alpha,S_n}^{\circ} \right \|_{\ell^\infty\to \ell^\infty}=(n-1)^\alpha.$$	
\end{theorem}

\begin{proof}
Let $a\in V_{S_n}$ be the center of $S_n$.  Then we have
\[
M_{\alpha,S_n}^{\circ}f(a)=\max_{r=0,1}\frac{1}{|S(a,r)|^{1-\alpha}}\sum_{x\in S(a,r)}|f(x)|
\]
and
\[
M_{\alpha,S_n}^{\circ}f(v)=\max_{r=0,1,2}\frac{1}{|S(v,r)|^{1-\alpha}}\sum_{x\in S(v,r)}|f(x)|
\]
for any vertex $v\in V_{S_n}\setminus\{a\}$.

Choose $f=\delta_a$, then for $0<q<\infty$,
\begin{align*}
	\left\|M_{\alpha,S_n}^{\circ}\delta_a\right\|_q
	=\bigg((M_{\alpha,S_n}^{\circ}\delta_a(a))^q+\sum_{v\in V_{S_n}\setminus\{a\}}(M_{\alpha,S_n}^{\circ}\delta_a(v))^q\bigg)^{1/q}.
\end{align*}

A simple calculation shows that $ M_{\alpha,S_n}^{\circ}\delta_a (a)= 1$ and
\[
M_{\alpha,S_n}^{\circ}\delta_a (v)= \max_{r=0,1,2}\frac{1}{|S(v,r)|^{1-\alpha}}\sum_{x\in S(v,r)}\delta_a(x)=\frac{1}{|S(v,1)|^{1-\alpha}}\sum_{x\in S(v,1)}\delta_a(x)=1
\]
for $ v\in V_{S_n}\setminus\{a\}$, and thus
\begin{align*}
	\left \| M_{\alpha,S_n}^{\circ}\delta_a\right \|_{q}=n^{1/q}.
\end{align*}
Thus we obtain the lower bound as follows:
\[
\left \| M_{\alpha,S_n}^{\circ}\right \|_{\ell^p\to \ell^q}\geq \frac{\left \| M_{\alpha,S_n}^{\circ}\delta_a\right \|_{q}}{\left \|\delta_a\right \|_{p}} =n^{1/q},
\ 0<q<\infty.
\]

By Theorem $\ref{thm2-3}$(i)-(iii), the identity $\Vert M_{\alpha,S_n}^{\circ}\Vert_{\ell^p \to \ell^q} = n^{1/q}$ holds for the cases: (a) $0 \leq \alpha < 1$, $0 < p \leq q \leq 1$; (b) $0< \alpha < 1$, $1 < p \leq q < \infty$ with $p \leq 1 / \alpha$; (c) $\alpha = 0$, $1 < p \leq q < \infty$; (d) $0 \leq \alpha < 1$, $0 < p \leq 1 < q < \infty$.

If $1<p\leq q<\infty$ and $p> {1}/{\alpha}$, then we have
\begin{align*}
	M_{\alpha,S_n}^{\circ}f(a)&=\max\Big\{|f(a)|,(n-1)^{\alpha-1}\sum_{x\in V_{S_n}\setminus\{a\}}|f(x)|\Big\}\\
	&\leq \max\big\{|f(a)|,(n-1)^{\alpha-{1/p}}\left\|f\right\|_p\big\}\\
     &= (n-1)^{\alpha-{1/p}}\left\|f\right\|_p,
\end{align*}
and for $v\in V_{S_n}\setminus\{a\}$, we have
{\allowdisplaybreaks\begin{align*}
	M_{\alpha,S_n}^{\circ}f(v)&=\max\Big\{|f(v)|,|f(a)|,(n-2)^{\alpha-1}\sum_{x\in V_{S_n}\setminus\{v,a\}}|f(x)|\Big\}\\
	&\leq \max\big\{|f(v)|,|f(a)|,(n-2)^{\alpha-{1/p}}\left\|f\right\|_p\big\}\\
&= (n-2)^{\alpha-{1/p}}\left\|f\right\|_p.
\end{align*}}

Assume that $\left\|f\right\|_p=1$. For $p>{1}/{\alpha}$ and  $x\in V_{S_n}$, we have
\[
\left \| M_{\alpha,S_n}^{\circ} \right \|_{\ell^p\to \ell^q}\leq \big((n-1)^{q(\alpha-{1/p})}+(n-1)(n-2)^{q(\alpha-{1/p})}\big)^{1/q}.
\]

For $0<\alpha<1$, by a similar method as that of proof (iv) of Theorem $\ref{thm2-4}$, we have
\begin{align}\label{equ2-6}
	\left\|M_{\alpha,S_n}^{\circ}\right\|_{\ell^p\to \ell^\infty}\leq \left\{\begin{matrix}
		1 & 0<p\leq1/\alpha,\\
		(n-1)^{\alpha-{1/p}} & 1/\alpha\leq p<\infty.
	\end{matrix}\right.
\end{align}

Taking two test functions $f=\delta_{i}$ for $i\neq a$ and $g=1-\delta_a$, we have
\begin{align*}
	\left\|M_{\alpha,S_n}^{\circ}\right\|_{\ell^p\to \ell^\infty}\geq \max\{1,(n-1)^{\alpha-{1/p}}\},
\end{align*}
which means that ${\rm(iv)}$ holds.

The proof of part (v) is omitted, as it mirrors the argument for Theorem $\ref{thm2-4}$(v).
\end{proof}

An application of the strong-type estimate for $M^{\circ}_{\alpha,S_n}$, presented above, yields the optimal norm of the operator in the weak-type setting. It is notationally convenient to set $1/0=\infty$.
\begin{corollary}\label{cor2-6}
Let $n\geq 3$, $0\leq \alpha<1$, and suppose $0<p\leq1/\alpha$,  $0< q\leq \infty$, and $p\leq q$. Then
$$\left\|M_{\alpha,S_n}^{\circ}\right\|_{\ell^{p}\to \ell^{q,\infty}}=n^{1/q}.$$	
\end{corollary}

\begin{proof}
Let the vertex $a$ be the center of $S_n$, and take $f=\delta_a$, then for all $v\in V_{S_n}$, $M_{\alpha,S_n}^{\circ}\delta_a(v)=1$ 
and
$$\{v\in V_{S_n}:M_{\alpha,S_n}^{\circ}\delta_a(v)>\lambda\}=\left\{\begin{matrix}
	V_{S_n} &\lambda\in (0,1), \\[0.1cm]
	\emptyset &\lambda\in [1,\infty).
\end{matrix}\right.$$
Thus we have
\begin{align*}
	\left\|M_{\alpha,S_n}^{\circ}(\delta_a)\right\|_{\ell^{q,\infty}}&=\sup_{\lambda>0}\lambda|\{v\in V_{S_n}:\ M_{\alpha,S_n}^{\circ}\delta_a(v)>\lambda\}|^{1/q}=n^{1/q},
\end{align*}	
which implies that
\[\left\|M_{\alpha,S_n}^{\circ}\right\|_{\ell^p\to \ell^{q,\infty}}\geq n^{1/q}.\]

Using Theorem \ref{thm2-6}, the upper bound can be obtained by
\[
\left\|M_{\alpha,S_n}^{\circ}\right\|_{\ell^{p}\to \ell^{q,\infty}}\leq \left\|M_{\alpha,S_n}^{\circ}\right\|_{\ell^p\to \ell^q}=n^{1/q}.
\]	
\end{proof}	

\begin{remark}
If $\alpha = 0$ and $p=q$, then Theorems $\ref{thm2-3}$, $\ref{thm2-4}$, $\ref{thm2-6}$ and Corollary $\ref{cor2-6}$ reduce to the results in \cite{LZ20212}.	
\end{remark}

\section{Characterizations  of extremal graphs}\label{s3}
This section focuses on characterizing extremal graphs via the operator norm $\left \| M_{\alpha,G}^{\circ} \right \|_{\ell^p\to \ell^q}$. We present three main theorems that precisely identify the graphs that maximize the operator norm under different constraints on the parameters $\alpha$, $p$, and $q$. Together, these results reveal a profound interplay between the graph structure and the norm behavior of the spherical maximal operator, as seen in the crucial role of high-degree vertices.

The following result establishes a characterization of the extremal graphs for which the operator norm $\left \| M_{\alpha,G}^{\circ} \right \|_{\ell^p\to \ell^q}$ attains its lower bound.
\begin{theorem}
Let $G=(V_G, E_G)$ be a finite connected graph with $|V_G|=n\geq 3$. Let $0\leq \alpha<1$, $0<p\leq 1$, $0<q<\infty$ and $p\leq q$. Then  $\left \| M_{\alpha,G}^{\circ} \right \|_{\ell^p\to \ell^q}=(1+(n-1)^{q(\alpha-1)+1})^{{1}/{q}}$ if and only if $G=K_n$.	
\end{theorem}

\begin{proof}
The sufficiency follows from Theorem \ref{thm2-4}(i) and (iii).

We now prove the necessity. To reach a contradiction, we suppose $G\ne K_n$. Then there exist two vertices $u,v\in V_G$ such that $d_G(u,v)>1$.  Set
$J=S(u,1)$ and $K=S(v,1)$. We consider the following two cases.

{\rm(i)} $\min\{|J|,|K|\}\leq \frac{n-2}{2}$.
We may assume $|J| \leq \frac{n-2}{2}$. By the connectivity of $G$, there exists a vertex $j \in J$ such that $d_G(j, u) = 1$. Now take $f = \delta_j$. Then for $0\leq \alpha<1$, we have
\begin{align*}
	\left\|M_{\alpha,G}^{\circ}\delta_j\right\|_q
	&=\bigg((M_{\alpha,G}^{\circ}\delta_j(j))^q+(M_{\alpha,G}^{\circ}\delta_j(u))^q+\sum_{x\in V_{G}\setminus\{j,u\}}(M_{\alpha,G}^{\circ}\delta_j(x))^q\bigg)^{{1}/{q}}\\
	&\geq \Big(1+|J|^{q(\alpha-1)}+(n-2)(n-1)^{q(\alpha-1)}\Big)^{{1}/{q}}\\
	&\geq \Big(1+2^{q(1-\alpha)}(n-2)^{q(\alpha-1)}+(n-2)(n-1)^{q(\alpha-1)}\Big)^{{1}/{q}}\\
	&>\left(1+(n-1)^{q(\alpha-1)+1}\right)^{{1}/{q}}.
\end{align*}

{\rm(ii)} $\min\{|J|,|K|\}> \frac{n-2}{2}$.  This means that $J\bigcap K\ne \emptyset$. Take $j \in J \bigcap K$ and let $f = \delta_j$. Then, for $0 \leq \alpha < 1$, it holds that
\begin{align*}
\left\|M_{\alpha,G}^{\circ}\delta_j\right\|_q 
&=\bigg(\sum_{x\in \{j,u,v\}}(M_{\alpha,G}^{\circ}\delta_j(x))^q+\sum_{x\in V_{G}\setminus\{j,u,v\}}(M_{\alpha, G}^{\circ}\delta_j(x))^q\bigg)^{{1}/{q}}\\
&\geq \left(1+|J|^{q(\alpha-1)}+|K|^{q(\alpha-1)}+(n-3)(n-1)^{q(\alpha-1)}\right)^{{1}/{q}}.
\end{align*}

Since  $|J|\leq n-2$ and $|K|\leq n-2$, then for  $j\in J\bigcap K$, we have
\begin{align*}
	\left\|M_{\alpha,G}^{\circ}\delta_j\right\|_q&\geq \left(1+2(n-2)^{q(\alpha-1)}+(n-3)(n-1)^{q(\alpha-1)}\right)^{{1}/{q}}\\
	&>\left(1+(n-1)^{q(\alpha-1)+1}\right)^{{1}/{q}}.
\end{align*}

Combining these two cases, we have
\begin{align*}
	\left \| M_{\alpha,G}^{\circ} \right \|_{\ell^p\to \ell^q}&\geq \left \| M_{\alpha,G}^{\circ}\delta_j \right \|_{q}>(1+(n-1)^{q(\alpha-1)+1})^{{1}/{q}},
\end{align*}
which yields a contradiction.
\end{proof}

The following result establishes a characterization of the extremal graphs for which the operator norm $\left \| M_{\alpha,G}^{\circ} \right \|_{\ell^p\to \ell^q}$ attains its upper bound. We prove that the equality $\left \| M_{\alpha,G}^{\circ} \right \|_{\ell^p\to \ell^q} = n^{1/q}$ occurs if and only if the graph is isomorphic to the star graph $S_n$, revealing an intimate connection between the analytic properties of the maximal operator and the geometric structure of the underlying graph.

\begin{theorem}\label{thm3-2}
Let $G = (V_G, E_G)$ be a finite connected graph with $|V_G| = n \geq 3$. Then
\[
\left\| M_{\alpha,G}^{\circ} \right\|_{\ell^p \to \ell^q} = n^{1/q}
\]
if and only if $G \cong S_n$, under the following conditions:

{\rm (i)} For $0 < \alpha < 1$, where $0 < p \leq 1/\alpha$, $0 < q < \infty$, and $q \geq p$;

{\rm (ii)} For $\alpha = 0$, where $0 < p \leq q < \infty$.
\end{theorem}

To establish this characterization, we first develop some technical tools. The following lemma provides precise pointwise estimates for the maximal operator $M_{\alpha,G}^{\circ}$, which will play a crucial role in the proof of Theorem \ref{thm3-2} by allowing us to analyze the extremal functions that achieve the operator norm.

\begin{lemma}\label{lem3-3}
 Let $G = (V_G, E_G)$ be a finite connected graph with $|V_G| = n\geq 3$, and let $f: V_G \to \mathbb{R}$ satisfy $\lVert f \rVert_p = 1$. Then

{\rm (i)} For $0 < \alpha < 1$, $0 < p \leq 1/\alpha$, $0 < q < \infty$ and $q \geq p$, we have
\[\lVert M_{\alpha,G}^{\circ} f \rVert_{q} = n^{1/q}\] if and only if
\[
M_{\alpha,G}^{\circ} f(v) = 1 \quad \text{for all } v \in V_G.
\]

For $\alpha = 0$ and $0 < p \leq q < \infty$, we have \[\lVert M_{G}^{\circ} f \rVert_{q} = n^{1/q}\] if and only if
\[
M_{G}^{\circ} f(v) = 1 \quad \text{for all } v \in V_G.
\]

{\rm (ii)} For $0 < \alpha < 1$, $0 < p \leq 1/\alpha$, $0 < q < \infty$ with $q \geq p$, the following holds: for every vertex $v \in V_G$, there exists some $r_v \in \{0, 1, \ldots, n-1\}$ such that
\begin{align}\label{max}
M_{\alpha,G}^{\circ} f(v) = |S(v, r_v)|^{\alpha - 1} \sum_{x \in S(v, r_v)} |f(x)|.
\end{align}
Moreover, the equality $M_{\alpha,G}^{\circ} f(v) = 1$ holds if and only if one of the following conditions is satisfied:
\begin{itemize}
\item[(a)]  When $0 < p < 1/\alpha$, it holds that $|S(v, r_v)| = 1$, and
\begin{align}\label{fff}
	|f(x)| =
	\begin{cases}
		1 & \text{if } x \in S(v, r_v), \\
		0 & \text{if } x \in V_G\setminus S(v, r_v).
	\end{cases}
\end{align}
\item[(b)] When $p = 1/\alpha$, we have
\begin{align}\label{ff3}
	|f(x)| =
	\begin{cases}
		|S(v, r_v)|^{-1/p} & \text{if } x \in S(v, r_v), \\
		0 & \text{if } x \in V_G\setminus S(v, r_v).
	\end{cases}
\end{align}	
\end{itemize}

{\rm (iii)} For $\alpha = 0$ and $0 < p \leq q < \infty$, the equality $M_{G}^{\circ} f(v) = 1$ holds if and only if
$
|S(v, r_v)| = 1
$
and $f$ satisfies \eqref{fff}.
\end{lemma}

\begin{proof}[Proof of Lemma  \ref{lem3-3}]
The proof of part (i) is straightforward by applying \eqref{equ2-1}, \eqref{p<1} for $0<p\leq 1$, and \eqref{p>1} for $1<p\leq 1/\alpha$ when $0<\alpha<1$, and by using \eqref{equ2-1}, \eqref{p<1} for  $0<p\leq 1$, and \eqref{equ77-2} for $1<p<\infty$ when $\alpha=0$.

We now proceed to part (ii). For $0 < \alpha < 1$, it follows from \eqref{equ2-1} that for any vertex $v\in V_{G}$, there exists some $r_v \in \{0, 1, \ldots, n-1\}$ such that
\[
M_{\alpha,G}^{\circ} f(v) = |S(v, r_v)|^{\alpha - 1} \sum_{x \in S(v, r_v)} |f(x)|.
\]

For  $0 < p \leq 1$, by \eqref{p<1}, we have that $M_{\alpha,G}^{\circ}f(v) = 1$ if and only if
\begin{align*}
|S(v, r_v)| = 1\ \text{and}\ \sum_{x \in S(v, r_v)} |f(x)|^{p}=1.
\end{align*}

For $1 < p \leq 1/\alpha$, it follows from \eqref{p>1} that
$M_{\alpha,G}^{\circ}f(v) = 1$ if and only if
\begin{align*}
|S(v, r_v)|^{\alpha - 1/p} = 1 \quad \text{and} \quad \sum_{x \in S(v, r_v)} |f(x)|^{p} = 1.
\end{align*}
This implies that for $1 < p < 1/\alpha$, we have
\begin{align*}
|S(v, r_v)| = 1 \quad \text{and} \quad \sum_{x \in S(v, r_v)} |f(x)|^{p} = 1,
\end{align*}
while for $p = 1/\alpha$, we obtain
\begin{align}\label{gg}
\sum_{x \in S(v, r_v)} |f(x)|^{p} = 1.
\end{align}

For $0 < p < 1/\alpha$, if we assume that $S(v, r_v) = \{t\}$, then $t$ is the unique vertex satisfying \eqref{fff}.

Now consider the case $p = 1/\alpha$. By \eqref{p>1} and \eqref{gg}, equality holds in
\[
M_{\alpha,G}^{\circ}f(v)\leq |S(v,r_v)|^{\alpha-{1/p}}\bigg(\sum_{x\in S(v,r_v)} |f(x)|^p\bigg)^{1/p}
\]
if and only if $|f(x)|^p =|S(v, r_v)|^{-1}$ and then condition $(b)$ is satisfied.

The argument for part (iii) follows in a completely analogous manner to that of part (ii), noting that in this case the value $p = \infty$ is excluded when $\alpha = 0$.
\end{proof}

\begin{proof} [Proof of Theorem  \ref{thm3-2}]
The sufficiency follows from conditions (i), (ii), and (iii) of Theorem \ref{thm2-6}.

The necessity will be proved by contradiction. Consider a function $f: V_G \to \mathbb{R}_+$ such that $\|f\|_p = 1$. Assume that $G$ is a finite connected graph and the following norm equality holds: if $0 < \alpha < 1$, with parameters satisfying $0 < p \leq 1/\alpha$, $0 < q < \infty$, and $q \geq p$,
\[
\left\| M_{\alpha,G}^{\circ} \right\|_{\ell^p \to \ell^q} = n^{1/q};
\]
if $\alpha = 0$, with parameters $0 < p \leq q < \infty$,
\[
\left\| M_{G}^{\circ} \right\|_{\ell^p \to \ell^q} = n^{1/q}.
\]
By Lemma \ref{lem3-3}(i), we have
\[\lVert M_{\alpha,G}^{\circ} f \rVert_{q} = n^{1/q}\] if and only if
\[
M_{\alpha,G}^{\circ} f(v) = 1 \quad \text{for all } v \in V_G.
\]
We now assume $G \not\cong S_n$. We first consider the case $0<\alpha<1$, which is analyzed in two subcases.

{\bf Case 1.} In this case, we denote by $V_G=\{u, v_1, v_2, \dots, v_{n-1}\}$. Assume that $\deg(u) = n - 1$ and $\deg(v_{k_0}) \geq 2$ for some $k_0 \in \{1, 2, \dots, n-1\}$. Then $d_{G}(u,v_i)=1$ for all $i=1,2,\ldots,n-1$.

Since $M_{\alpha,G}^{\circ} f(u) = 1 $, by Lemma \ref{lem3-3}(ii), there exists $r_{u}\in \{0,1\}$ such that
\begin{align*}
M_{\alpha,G}^{\circ} f(u) = |S(u, r_u)|^{\alpha - 1} \sum_{x \in S(u, r_u)} |f(x)|
\end{align*}
holds.

For $0<p<1/\alpha$, by Lemma \ref{lem3-3}(ii),  we have
$|S(u,r_u)|=1$ which implies that $r_u=0$ and $S(u,r_u)=\{u\}$.
Hence we have $|f(u)|=1$ and $f(v)=0$ for all $v\in V_G\setminus\{u\}$.

For $p=1/\alpha$, by Lemma \ref{lem3-3}(ii),  we have
\begin{align}\label{ff33}
|f(x)| =
\begin{cases}
|S(u, r_u)|^{-1/p} & \text{if } x \in S(u, r_u), \\
 0 & \text{if } x \in V_G\setminus S(u, r_u).
\end{cases}
\end{align}
If $r_u=0$, then \eqref{ff33} is
\begin{align*}
|f(x)| =
\begin{cases}
1 & \text{if } x=u, \\
 0 & \text{if } x \in V_G\setminus \{u\}.
\end{cases}
\end{align*}
If $r_u=1$, then \eqref{ff33} is
\begin{align*}
|f(x)| =
\begin{cases}
(n-1)^{-1/p} & \text{if } x \in V_G\setminus \{u\}, \\
 0 & \text{if } x=u.
\end{cases}
\end{align*}

We conclude that for $0<p\leq 1/\alpha$ and $r_u=0$,  we have
\begin{align*}
|f(x)| =
\begin{cases}
1 & \text{if } x=u, \\
 0 & \text{if } x \in V_G\setminus \{u\},
\end{cases}
\end{align*}
and for $p=1/\alpha$ and $r_u=1$, we have
\begin{align}\label{lsz-1}
|f(x)| =
\begin{cases}
(n-1)^{-1/p} & \text{if } x \in V_G\setminus \{u\}, \\
 0 & \text{if } x=u.
\end{cases}
\end{align}

It follows from  $M_{\alpha,G}^{\circ} f(v_{k_0}) = 1$ that there exists some $r_{v_{k_0}} \in \{0, 1, \ldots, n-1\}$ such that
\begin{align}\label{max22}
M_{\alpha,G}^{\circ} f(v_{k_0}) = |S(v_{k_0}, r_{v_{k_0}})|^{\alpha - 1} \sum_{x \in S(v_{k_0}, r_{v_{k_0}})} |f(x)|.
\end{align}

If $0<p\leq 1/\alpha$ and $r_u=0$,  then we have $r_{v_{k_0}}=1$ because  $d_{G}(u,v_{k_0})=1$. Using the hypothesis $\deg(v_{k_0}) \geq 2$, we have  for $0<\alpha< 1$,
\[
1=M_{\alpha,G}^{\circ}f(v_{k_0})=|S(v_{k_0},1)|^{\alpha-1}|f(u)|\leq 2^{\alpha-1}<1,
\]
which yields a contradiction.

We now turn to the case $p=1/\alpha$ and $r_u=1$. If $r_{v_{k_0}}=0$, then by \eqref{lsz-1}, we have
\[
M_{\alpha,G}^{\circ}f(v_{k_0})=|f(v_{k_0})|=(n-1)^{-1/p}<1.
\]
If $r_{v_{k_0}}\geq 1$, then $v_{k_0}\notin S(v_{k_0},r_{v_{k_0}})$. Since $|f(u)|=0$, we have
\begin{align*}
    M_{\alpha,G}^{\circ}f(v_{k_0})&=|S(v_{k_0},r_{v_{k_0}})|^{\alpha-1}\sum_{x\in S(v_{k_0},r_{v_{k_0}})}|f(x)|\\
    &\leq |S(v_{k_0},r_{v_{k_0}})|^{\alpha-1}(n-2)(n-1)^{-1/p}\\
    &\leq (n-1)^{\alpha-1}(n-2)(n-1)^{-1/p}\\
    &<1,
\end{align*}
which yields a contradiction.

{\bf Case 2.}
Assume that $\deg(v)\leq n-2$ for every vertex $v \in V_G$. Let $D = \max\{\deg(v) : v \in V_G\}$, and suppose there exists a vertex $u \in V_G$ such that $\deg(u) = D$. Since $G$ is connected, it follows that $D \geq 2$. By Lemma \ref{lem3-3}(ii), there exists $r_u \in \{0, 1, \dots, n-1\}$ such that
\begin{align*}
    M_{\alpha,G}^{\circ} f(u) = |S(u, r_u)|^{\alpha - 1} \sum_{x \in S(u, r_u)} |f(x)|.
\end{align*}
Note that $M_{\alpha,G}^{\circ} f(u) = 1$. Applying Lemma \ref{lem3-3}(ii) again, we obtain the following cases:

\begin{enumerate}
\item[(a)] If $r_u = 0$ and $0 < p \leq 1/\alpha$, then $|f(u)| = 1$ and $f(v) = 0$ for all $v \in V_G \setminus \{u\}$.

\item[(b)] If $r_u = 1$ and $0 < p < 1/\alpha$, then $|S(u, 1)| = 1$. However, this contradicts $\deg(u) \geq 2$, implying that $r_u \geq 2$ when $0 < p < 1/\alpha$.

\item[(c)] If $r_u \geq 2$ and $0 < p < 1/\alpha$, then $|S(u, r_u)| = 1$. Denote $S(u, r_u) = \{t_{r_u}\}$. Then $|f(t_{r_u})| = 1$ and $f(v) = 0$ for all $v \in V_G \setminus \{t_{r_u}\}$.

\item[(d)] If $r_u \geq 1$ and $p = 1/\alpha$, then we have
\begin{align*}
	|f(x)| =
	\begin{cases}
		|S(u, r_u)|^{-1/p} & \text{if } x \in S(u, r_u), \\
		0 & \text{otherwise},
	\end{cases}
\end{align*}
when $|S(u, r_u)| \geq 2$; and
\begin{align*}
	S(u, r_u) = \{t_{r_u}\} \text{ with } |f(t_{r_u})| = 1, \text{ and } f(v) = 0 \text{ for all } v \in V_G \setminus \{t_{r_u}\},
\end{align*}
when $|S(u, r_u)| = 1$.
\end{enumerate}

Therefore, we conclude the following:
\begin{enumerate}
\item For $0 < p \leq 1/\alpha$ with $r_u = 0$,
\begin{align*}
	|f(x)| =
	\begin{cases}
		1 & \text{if } x = u, \\
		0 & \text{if } x \in V_G \setminus \{u\}.
	\end{cases}
\end{align*}

\item For $0 < p \leq 1/\alpha$ with $r_u \geq 2$ and $|S(u, r_u)| = 1$,
\begin{align*}
	|f(x)| =
	\begin{cases}
		1 & \text{if } x \in S(u, r_u), \\
		0 & \text{if } x \in V_G \setminus S(u, r_u).
	\end{cases}
\end{align*}

\item For $p = 1/\alpha$ with $r_u \geq 1$ and $|S(u, r_u)| \geq 2$,
\begin{align*}
	|f(x)| =
	\begin{cases}
		|S(u, r_u)|^{-1/p} & \text{if } x \in S(u, r_u), \\
		0 & \text{if } x \in V_G \setminus S(u, r_u).
	\end{cases}
\end{align*}
\end{enumerate}

We now consider the vertex set $V_G = \{u, v_1, v_2, \ldots, v_{n-1}\}$. Let $N_G(u)$ denote the neighborhood of vertex $u$. For notational convenience, we define $N_G(u) = \{v_i \in V_G : 1 \leq i \leq D\}$ and $V_G \setminus (N_G(u)\cup \{u\}) = \{v_i \in V_G : D+1 \leq i \leq n-1\}$. By the connectivity of $G$, there exist two vertices $v_{i_0} \in N_G(u)$ and $v_{j_0} \in V_G \setminus  (N_G(u)\cup \{u\})$ such that $v_{j_0} \sim v_{i_0}$. Observe that both $u$ and $v_{j_0}$ belong to $N_G(v_{i_0})$, which implies $|S(v_{i_0}, 1)| \geq 2$.

Given that $0 < \alpha < 1$, when $0 < p \leq 1/\alpha$ and $r_u = 0$, we obtain
\[
M_{\alpha,G}^{\circ} f(v_{i_0}) = |S(v_{i_0}, 1)|^{\alpha - 1} |f(u)| \leq 2^{\alpha - 1} < 1.
\]

For the case $0 < p \leq 1/\alpha$ with $r_u \geq 2$ and $|S(u, r_u)| = 1$, let $S(u, r_u) = \{t_{r_u}\}$, and let $P$ be a shortest path from $u$ to $t_{r_u}$. Denote by $y$ the vertex adjacent to $t_{r_u}$ on path $P$. Since $|S(y, 1)| \geq 2$, we have
\[
M_{\alpha,G}^{\circ} f(y) = |S(y, 1)|^{\alpha - 1} |f(t_{r_u})| \leq 2^{\alpha - 1} < 1.
\]

Now consider the case $p = 1/\alpha$ with $r_u \geq 1$ and $|S(u, r_u)| \geq 2$. Select a vertex $w\in S(u, r_u)$. Then there exists $r_w \in \{0, 1, \ldots, n-1\}$ such that
\[
M_{\alpha,G}^{\circ} f(w) = |S(w, r_w)|^{\alpha - 1} \sum_{x \in S(w, r_w)} |f(x)|.
\]
We examine two subcases:

(1) If $r_w = 0$, then $S(w, r_w) = \{w\}$ and $|f(w)| = |S(u, r_u)|^{-1/p}$, yielding
    \[
    M_{\alpha,G}^{\circ} f(w) = |S(w, r_w)|^{\alpha - 1} |f(w)| = |S(u, r_u)|^{-\alpha} < 1.
    \]

(2) If $r_w \geq 1$, then from the previous analysis we have $|f(v)| = |S(u, r_u)|^{-1/p}$ for $v \in S(u, r_u)$ and $f(v) = 0$ for $v \in V_G \setminus S(u, r_u)$. We further consider the following two cases.

If $|S(w, r_w)| \geq |S(u, r_u)|$, noting that $w \in S(u, r_u)$ but $w \notin S(w, r_w)$ (since $r_w \geq 1$), we obtain
\begin{align*}
	M_{\alpha,G}^{\circ} f(w) &= |S(w, r_w)|^{\alpha - 1} \sum_{x \in S(w, r_w)} |f(x)| \\
	&= |S(w, r_w)|^{\alpha - 1} \sum_{x \in S(w, r_w) \cap S(u, r_u)} |f(x)| \\
	&\leq |S(w, r_w)|^{\alpha - 1} (|S(u, r_u)| - 1) |S(u, r_u)|^{-1/p} \\
	&\leq |S(u, r_u)|^{\alpha - 1} (|S(u, r_u)| - 1) |S(u, r_u)|^{-1/p} \\
	&< 1.
\end{align*}

If $|S(w, r_w)| < |S(u, r_u)|$, then
\[
M_{\alpha,G}^{\circ} f(w) \leq |S(w, r_w)|^{\alpha - 1} |S(w, r_w)| |S(u, r_u)|^{-1/p} < 1.
\]

In all cases, we arrive at a contradiction to the hypothesis that $M_{\alpha, G}^{\circ} f(v) = 1$ for all $v \in V_G$.

The case $\alpha = 0$ can be handled by an argument analogous to that for the case $0 < \alpha < 1$ (with the noted exclusion of $p = \infty$). The details are a straightforward adaptation and are therefore omitted.
\end{proof}

The following result provides a complete characterization of the extremal graphs for which the operator norm $\left\|M_{\alpha,G}^{\circ}\right\|_{\ell^p\to \ell^\infty}$ attains its maximum possible value. It establishes that the equality $\left\|M_{\alpha,G}^{\circ}\right\|_{\ell^p\to \ell^\infty}=(n-1)^{\alpha-1/p}$ holds if and only if the graph contains a universal vertex, explicitly connecting the analytic operator norm with the combinatorial structure of the graph.

\begin{theorem}\label{thm3-4}
Let $G=(V_G,E_G)$ be a finite connected graph with $|V_G|=n\geq 3$. Let $0\leq \alpha<1$ and $1/\alpha\leq p\leq\infty$. Then $\left\|M_{\alpha,G}^{\circ}\right\|_{\ell^p\to \ell^\infty}=(n-1)^{\alpha-1/p}$ if and only if there exists a vertex $u\in V_G$ such that $\deg(u)=n-1$.
\end{theorem}

\begin{proof}
We begin by proving the sufficiency. Assume that there exists a vertex $u \in V_G$ with $\deg(u) = n - 1$. Choose the function $f = 1 - \delta_u$. Then we have $\|f\|_p = (n - 1)^{1/p}$.

For any vertex $v \in V_G \setminus \{u\}$, there exists some $r_v \in \{0, 1, \dots, n-1\}$ such that
\begin{align*}
    M^{\circ}_{\alpha,G} f(v) &= \frac{1}{|S(v, r_v)|^{1 - \alpha}} \sum_{x \in S(v, r_v)} |f(x)| \\
    &=|S(v, r_v)|^{\alpha} \\
    &\leq (n - 1)^{\alpha}.
\end{align*}

Since $\deg(u) = n - 1$, we have
\begin{align*}
    M^{\circ}_{\alpha,G} f(u) &= \max_{r = 0, 1} \frac{1}{|S(u, r)|^{1 - \alpha}} \sum_{x \in S(u, r)} |f(x)| \\
    &= \max \bigg\{ |f(u)|,\ (n - 1)^{\alpha - 1} \sum_{x \in S(u, 1)} |f(x)| \bigg\} \\
    &= (n - 1)^{\alpha}.
\end{align*}

Combining these estimates, we obtain for $f = 1 - \delta_u$,
\[
\left\| M_{\alpha,G}^{\circ} f \right\|_{\ell^\infty} = (n - 1)^{\alpha},
\]
and therefore
\[
\left\| M_{\alpha,G}^{\circ} \right\|_{\ell^p \to \ell^\infty} \geq \frac{(n - 1)^{\alpha}}{(n - 1)^{1/p}} = (n - 1)^{\alpha - 1/p}.
\]

By Theorem \ref{thm2-3}(iv), we also have the upper bound
\[
\left\| M_{\alpha,G}^{\circ} \right\|_{\ell^p \to \ell^\infty} \leq (n - 1)^{\alpha - 1/p}.
\]
Hence, we conclude that
\[
\left\| M_{\alpha,G}^{\circ} \right\|_{\ell^p \to \ell^\infty} = (n - 1)^{\alpha - 1/p},
\]
which completes the proof of sufficiency.

Next, we prove the necessity. Let $\left\| M^{\circ}_{\alpha,G} \right\|_{\ell^p \to \ell^\infty} = (n - 1)^{\alpha - 1/p}$. Assume that $\deg(v) \leq n - 2$ for every vertex $v \in V_G$. We first observe that for any vertex $v$,
\begin{align*}
|S(v, 0)| = 1 \quad \text{and} \quad |S(v, 1)| = \deg(v) \leq n - 2.
\end{align*}

Since $G$ is connected, for any vertex $v$ there exists a vertex $v_0 \in N_G(v)$. Consequently, for any $r \geq 2$, we have
\[
|S(v, r)| \leq n - 2.
\]

Now, for any function $f \in \ell^p(V_G)$ and any vertex $v \in V_G$, there exists some $r_v \in \{0, 1, \dots, n-1\}$ such that
\begin{align*}
M^{\circ}_{\alpha,G} f(v) = \frac{1}{|S(v, r_v)|^{1 - \alpha}} \sum_{x \in S(v, r_v)} |f(x)|.
\end{align*}

We proceed by considering two cases based on the cardinality of $S(v, r_v)$:

\textbf{Case 1.} $|S(v, r_v)| = 1$.

Let $S(v, r_v) = \{t\}$. It is easy to see that
\begin{align*}
	M^{\circ}_{\alpha,G} f(v) = |f(t)| \leq \|f\|_p.
\end{align*}

\textbf{Case 2.} $2 \leq |S(v, r_v)| \leq n - 2$.

The H\"{o}lder inequality yields that for $\alpha \geq 1/p$, we have
\begin{align*}
M^{\circ}_{\alpha,G} f(v) \leq |S(v, r_v)|^{\alpha - 1/p} \|f\|_p  \leq (n - 2)^{\alpha - 1/p} \|f\|_p.	
\end{align*}

Combining both cases, we obtain the uniform estimate
\[
M^{\circ}_{\alpha,G} f(v) \leq \max\left\{1, (n - 2)^{\alpha - 1/p}\right\} \|f\|_p = (n - 2)^{\alpha - 1/p} \|f\|_p
\]
for all $v \in V_G$. This yields that the estimate $\left\| M^{\circ}_{\alpha,G} \right\|_{\ell^p \to \ell^\infty} \leq (n - 2)^{\alpha - 1/p}$, which contradicts the hypothesis that $\left\| M^{\circ}_{\alpha,G} \right\|_{\ell^p \to \ell^\infty} = (n - 1)^{\alpha - 1/p}$, thus completing the proof of necessity.
\end{proof}

\begin{remark}
The existence of a vertex $u \in V_G$ with $\deg(u) = n-1$ is a property shared by several graph families on $n$ vertices. These include, but are not limited to, complete graphs, star graphs, windmill graphs, and wheel graphs.
\end{remark}

\section{Optimal estimates for  $\left\Vert M_{P_n}^{\circ} \right\Vert_{\ell^p\to \ell^p}$} \label{s4}

In this section, we investigate the path $P_n$, which represents a simple yet fundamental non-symmetric structure. Unlike the previously studied complete and star graphs, the path exhibits linear connectivity and thus requires a finer analysis that depends heavily on the parity of $n$ and the location of vertices relative to the center, making it challenging to establish sharp bounds even for the spherical maximal operator on paths. Given this difficulty, we do not pursue the fractional case. The following theorem provides the exact values of the operator norm of $M_{P_n}^{\circ}$ for $0 < p \leq 1$ and $p = \infty$, and sharper two-sided bounds for $1 < p < \infty$.

\begin{theorem}\label{thm4-1}
Let $P_n = (V_{P_n}, E_{P_n})$ be the path with $n \geq 3$, where $V_{P_n} = \{1, 2, \ldots, n\}$,
$E_{P_n} = \{1 \sim 2, 2 \sim 3, \ldots, n-1 \sim n\}$, and define $m \geq 2$ such that:
$n = 2m-1$ if $n$ is odd, and $n = 2m$ if $n$ is even.

{\rm(i)} If $0 < p \leq 1$, then
\[
\left\Vert M_{P_n}^{\circ} \right\Vert_{\ell^p\to \ell^p} =
\begin{cases}
	\left( m + 1 + \dfrac{m - 2}{2^p} \right)^{1/p} & n = 2m - 1, \\[0.3cm]
	\left( m + 1 + \dfrac{m - 1}{2^p} \right)^{1/p} & n = 2m.
\end{cases}
\]

{\rm(ii)} If $1 < p < \infty$, then
\begin{align*}
	\left( m + 1 + \dfrac{m - 2}{2^p} \right)^{1/p} \leq \left\| M_{P_n}^{\circ} \right\|_{\ell^p\to \ell^p}\leq \left( \dfrac{5}{4} \left( m + 1 + \dfrac{m - 2}{2^p} \right) \right)^{1/p}
\end{align*}
for $n=2m-1$, and
\begin{align*}
	\left( m + 1 + \dfrac{m - 1}{2^p} \right)^{1/p} \leq \left\| M_{P_n}^{\circ} \right\|_{\ell^p\to \ell^p}\leq \left( \dfrac{5}{4} \left( m + 1 + \dfrac{m - 1}{2^p} \right) \right)^{1/p}
\end{align*}
for $n=2m$.

{\rm(iii)} If $p=\infty$, then
\begin{align*}
	\left\| M_{P_n}^{\circ} \right\|_{\ell^\infty\to \ell^\infty}=1.
\end{align*}

{\rm(iv)} If $0 < p < \infty$, then for $n=3$, we have
\[
\left\| M_{P_3}^{\circ} \right\|_{\ell^p\to \ell^p} = 3^{1/p}.
\]
\end{theorem}

We first give a direct calculation.

\begin{lemma}\label{lem4-2}
Let $P_n = (V_{P_n}, E_{P_n})$ be the path with $n \geq 3$, where $V_{P_n} = \{1, 2, \ldots, n\}$ and $E_{P_n} = \{1 \sim 2, 2 \sim 3, \ldots, n-1 \sim n\}$. Let $f: V_{P_n} \to \mathbb{R}$ and define its extension $\tilde{f}: \mathbb{Z}\to \mathbb{R}$  as
\begin{align}\label{equ4-1}
\tilde{f}(x) = \begin{cases}
	f(x) & \text{if } x\in V_{P_n},\\
	0    & \text{otherwise}.
\end{cases}	
\end{align}
Then it holds that
\[
M^{\circ}_{P_n}f(x) = \max_{r=0,1,\ldots,n-1} A_{P_{n},r}^{\circ}(|f|)(x), \ x\in V_{P_n},
\]
where the average function $A^{\circ}_{P_n,r}f$ of $f$ is defined by
\[
A_{P_{n},r}^{\circ}f(x)=\frac{1}{|S(x,r)|}\sum_{y\in S(x,r)} f(y),
\]
and
\begin{align}  \label{equ4-2}
A_{P_{n},r}^{\circ}(\vert f\vert) (x) =
\begin{cases}
	|f(x)| & \text{if } r=0,\\[5pt]
	\dfrac{|\tilde{f}(x-r)| + |\tilde{f}(x+r)|}{2} & \text{if } 0 < r \leq \min\lbrace d_{P_n}(x,1),\, d_{P_n}(x,n)\rbrace,\\[5pt]
	|\tilde{f}(x-r)| + |\tilde{f}(x+r)| & \text{if } \min\lbrace d_{P_n}(x,1),\, d_{P_n}(x,n)\rbrace < r \leq n-1.
\end{cases}	
\end{align}	
\end{lemma}	
	
\begin{lemma}\label{lem4-3} Let $1\leq p<\infty$. For any real number $\zeta\in \left[0, 1\right]$, let the function
\begin{align}\label{equ4-3}
	g(x)=x+2^{-p}\zeta\left(x^{1/p}+(1-x)^{1/p}\right)^p, \ x\in \left[1/2,1\right].
\end{align}
Then for $x\in \left[1/2,1\right]$, we have
$$g(x)\leq \frac{5}{4} g(1).$$	
\end{lemma}	
	
\begin{proof}[Proof of Lemma \ref{lem4-3}]
The case $\zeta = 0$ is trivial. Henceforth, we assume $0 < \zeta \leq 1$. For $1 \leq  p < \infty$, define $\tau = \lfloor p \rfloor$. Then $\tau \in \mathbb{N}_{+}$ and $1 \leq \tau < \infty$. By the power mean inequality, for $0 < \beta < \alpha$ and nonnegative real numbers $c$ and $d$, we have
\[
\bigg( \frac{c^\beta + d^\beta}{2} \bigg)^{1/\beta} \leq \bigg( \frac{c^\alpha + d^\alpha}{2} \bigg)^{1/\alpha}.
\]
Thus, for $1/2 \leq  x \leq 1$,
\[
g(x) \leq x + \zeta \left( \frac{1}{2} \left( x^{1/\tau} + (1-x)^{1/\tau} \right) \right)^\tau.
\]
On the other hand,
\[
g(1) = 1 + 2^{-p} \zeta  \geq 1 + 2^{-(\tau+1)}\zeta.
\]
Therefore, it suffices to prove that the inequality
\begin{align}\label{equ4-4}
	x + \zeta \bigg( \frac{1}{2} \bigg( x^{1/\tau} + (1-x)^{1/\tau} \bigg) \bigg)^\tau \leq \frac{5}{4} \big( 1 + 2^{-(\tau+1)}\zeta \big)
\end{align}
holds for all $0 < \zeta \leq 1$ and $x \in [1/2, 1]$.

To prove \eqref{equ4-4}, set $z = \left( \frac{1-x}{x} \right)^{{1}/{\tau}} \in [0,1]$. Then $x = \frac{1}{z^{\tau} + 1}$, and \eqref{equ4-4} becomes
\[
2^{\tau-2} (1 + 5z^{\tau}) + \zeta \left( \frac{5}{8} + \frac{5}{8} z^{\tau} - (1 + z)^{\tau} \right) \geq 0.
\]
Note that $0 < \zeta \leq 1$ and
\[
\frac{5}{8} + \frac{5}{8} z^{\tau} - (1 + z)^{\tau} \leq 1 + z^{\tau} - (1 + z)^{\tau} \leq 0.
\]
Therefore, to establish \eqref{equ4-4}, it suffices to prove that for all $z \in [0,1]$,
\begin{equation}\label{equ4-5}
	h(z) = 2^{\tau-2} (1 + 5z^{\tau}) + \frac{5}{8} + \frac{5}{8} z^{\tau} - (1 + z)^{\tau} \geq 0.
\end{equation}
Clearly, $h(0) > 0$.  For $0 < z \leq 1$, we have		
\[
h(z) =
\begin{cases}
\dfrac{1}{8} + \dfrac{17}{8}z, & \tau = 1, \\[0.3cm]
\dfrac{37}{8}z^2 - 2z + \dfrac{5}{8}, & \tau = 2, \\[0.3cm]
\dfrac{77}{8}z^3 - 3z^2 - 3z + \dfrac{13}{8}, & \tau = 3.
\end{cases}
\]
Obviously, $ h(z) > 0 $ for $0 < z \leq 1$.
		
For the case $\tau \geq 4$ and $0 < z \leq 1$, we proceed as follows.
A direct calculation leads to
\[
h'(z) = \tau z^{\tau-1} \left( 5 \cdot 2^{\tau-2} + \frac{5}{8} - \left(1 + \frac{1}{z}\right)^{\tau-1} \right).
\]
Define $u(z) = 5 \cdot 2^{\tau-2} + \frac{5}{8} - \left(1 + \frac{1}{z}\right)^{\tau-1}$. It is easy to verify that $u$ is increasing on $(0,1]$. Note that $\lim_{z \to 0} u(z) = -\infty < 0$ and $u(1) = \frac{5}{8} + 3 \cdot 2^{\tau-2} > 0$.
This implies that there is exactly one minimum point $z_0 \in (0,1)$ satisfying $h'(z_0) = 0$. Then,
\[
z_0 = {\bigg(\left( 5 \cdot 2^{\tau-2} + \frac{5}{8} \right)^{1/(\tau-1)} - 1\bigg)^{-1}},
\]
and the minimum value of $h(z)$ on $(0,1]$ is $h(z_0)$, which can be rewritten as
\[
h(z_0) = \frac{5}{8} + 2^{\tau-2} - \left( \frac{A}{A-1} \right)^{\tau-1},
\]
where $A = \left( 5 \cdot 2^{\tau-2} + \frac{5}{8} \right)^{1/(\tau-1)}$.

We can prove \eqref{equ4-5} using the fact that $2^{\tau-2} \geq \left( \frac{A}{A-1} \right)^{\tau-1}$ for $\tau \geq 4$. Specifically, it suffices to prove that for $\tau \geq 4$,
\begin{align}\label{equ4-6}
	5 \left( 2^{\frac{\tau-2}{\tau-1}} - 1 \right)^{\tau-1} \geq 1,
\end{align}
which is equivalent to
\[
\bigg( \frac{ \left( \frac{2}{5} \right)^{\frac{1}{\tau-1}} + 2^{\frac{1}{\tau-1}} }{2} \bigg)^{\tau-1} \leq 1.
\]
By the power mean inequality, it holds that
\[
\bigg( \frac{ \left( \frac{2}{5} \right)^{\frac{1}{\tau-1}} + 2^{\frac{1}{\tau-1}} }{2} \bigg)^{\tau-1} \leq \bigg( \frac{ \left( \frac{2}{5} \right)^{\frac{1}{3}} + 2^{\frac{1}{3}} }{2} \bigg)^{3} \leq 1.
\]
Therefore, we complete the proof of the lemma.		
\end{proof}

\begin{proof}[Proof of Theorem \ref{thm4-1}]
(i) For any $t\in V_{P_n}$, take $f=\delta_t$, then $M_{P_n}^{\circ}\delta_t(t)=1$. For any vertex $i\in V_{P_n}\setminus\{t\}$, according to Lemma \ref{lem4-2}, we calculate that
\[
d_{P_n}(i,1)=i-1,\qquad d_{P_n}(i,n)=n-i, \qquad d_{P_n}(t,i)=|t-i|.
\]
For $i\in V_{P_n}\setminus\{t\}$, we have
\begin{align*}
	A_{P_{n},r}^{\circ}(\delta_{t})(i)&=
	\left\{\begin{matrix}
		0	& r=0,\\
		(\tilde{\delta_{t}}(i-r)+\tilde{\delta_{t}}(i+r))/2	& 0<r\leq \min\{i-1,n-i\},\\
		\tilde{\delta_{t}}(i-r)+\tilde{\delta_{t}}(i+r)	& \min\{i-1,n-i\}<r\leq n-1,
	\end{matrix}\right.
\end{align*}
where $\tilde{\delta_{t}}$ is defined as \eqref{equ4-1}.
By Lemma \ref{lem4-2}, for $i\in V_{P_n}\setminus\{t\}$, we have
\begin{align*}
	M_{P_n}^{\circ}\delta_t(i)&=\max_{r=0,1,\ldots,n-1} A_{P_{n},r}^{\circ}\delta_{t}(i)
=A_{P_{n},d_{P_n}(t,i)}^{\circ}\delta_{t}(i)\\
&=\left\{\begin{matrix}
		\frac{1}{2}	& 1\leq d_{P_n}(t,i)\leq \min\{i-1,n-i\},\\
		1	& \min\{i-1,n-i\}<d_{P_n}(t,i)\leq n-1.
	\end{matrix}\right.
\end{align*}
To evaluate $M_{P_n}^{\circ}\delta_t(i)$, it is natural to distinguish between two cases: $i < t$ and $i > t$.

For $i<t$, we proceed by examining the two distinct cases: when $t$ is odd and when $t$ is even. If $t = 2k - 1$, where  $k \in \mathbb{N}_{+}$ and $2 \leq k \leq (n+1)/2$, then we have
\begin{align*}
	M_{P_n}^{\circ}\delta_t(i)=\begin{cases}
		\frac{1}{2}& i\geq k, \\
		1 &  i\leq k-1.
	\end{cases}
\end{align*}
If $t=2k$  with $k \in \mathbb{N}_{+}$ and $1 \leq  k \leq n/2$, then we have
\begin{align*}
	M_{P_n}^{\circ}\delta_t(i)=\begin{cases}
		\frac{1}{2}&  i\geq k+1, \\
		1 & i\leq k.
	\end{cases}
\end{align*}

For $i>t$, we consider the cases where $n+t$ is odd and even separately.
If $n+t=2l-1$ with $l\in \mathbb{N}_{+}$ and $2\leq l\leq n$, then \begin{align*}
	M_{P_n}^{\circ}\delta_t(i)=\begin{cases}
		\frac{1}{2} & i\leq  l-1, \\
		1 & i\geq l.
	\end{cases}
\end{align*}
If $n+t=2l$ with $l\in \mathbb{N}_{+}$ and $1\leq l\leq n-1$, then
\begin{align*}
	M_{P_n}^{\circ}\delta_t(i)=\begin{cases}
		\frac{1}{2} & i\leq  l, \\
		1 &  i\geq l+1.
	\end{cases}
\end{align*}

For $0<p<\infty$, to evaluate $\left\Vert M_{P_n}^{\circ}\delta_t \right \Vert_{p}$, we distinguish between the cases based on the cardinality $n$ of the vertex set $V_{P_n}$: when $n$ is odd and when $n$ is even.

{\bf Case 1.} $n = 2m-1$, where $m \in \mathbb{N}_{+}$ and $m \geq 2$.

If $t$ is odd, i.e., $t = 2k-1$ for some positive integer $1\leq k\leq m$, then $n+t$ is even. When $k=1$, based on the above estimates, we obtain that for $t=1$,
\begin{align*}
\left \| M_{P_n}^{\circ}\delta_1 \right \|_{p}&=\bigg((M_{P_n}^{\circ}\delta_{1}(1))^{p}+\sum_{2\leq i\leq m}(M_{P_n}^{\circ}\delta_{1}(i))^p+\sum_{m+1\leq i\leq 2m-1}(M_{P_n}^{\circ}\delta_{1}(i))^p\bigg)^{1/p}\\
&=\left(1+(m-1)2^{-p}+(m-1)\cdot 1\right)^{1/p}\\
&=\left(m+2^{-p}(m-1)\right)^{1/p}.	
\end{align*}
When $k=m$, we yield that for $t=2m-1$,
\begin{align*}
&\left \| M_{P_n}^{\circ}\delta_{2m-1} \right \|_{p}\\
&=\bigg(\sum_{1\leq i\leq m-1}(M_{P_n}^{\circ}\delta_{2m-1}(i))^p+\sum_{m\leq i\leq 2m-2}(M_{P_n}^{\circ}\delta_{2m-1}(i))^p+(M_{P_n}^{\circ}\delta_{2m-1}(2m-1))^{p}\bigg)^{1/p}\\
&=\left((m-1)\cdot 1+(m-1)2^{-p}+1\right)^{1/p}\\
&=\left(m+2^{-p}(m-1)\right)^{1/p}.	
\end{align*}
When $2\leq k\leq m-1$, we have
{\allowdisplaybreaks\begin{align*}
\left \| M_{P_n}^{\circ}\delta_t \right \|_{p}
&=\bigg(\sum_{1\leq i\leq k-1}(M_{P_n}^{\circ}\delta_{t}(i))^p+\sum_{k\leq i\leq 2k-2}(M_{P_n}^{\circ}\delta_{t}(i))^p+(M_{P_n}^{\circ}\delta_{t}(t))^{p}\\
&\quad +\sum_{2k\leq i\leq m+k-1}(M_{P_n}^{\circ}\delta_{t}(i))^p+\sum_{m+k\leq i\leq 2m-1}(M_{P_n}^{\circ}\delta_{t}(i))^p\bigg)^{1/p}\\
&=\left((k-1)\cdot 1+(k-1)2^{-p}+1+(m-k)2^{-p}+(m-k)\cdot 1\right)^{1/p}\\
&=\big(m+2^{-p}(m-1)\big)^{1/p}.		
\end{align*}}

Similarly, if $t$ is even, i.e., $t = 2k$ for some integer $k$ with $1 \leq k \leq m-1$, then $n + t$ is odd. When $k=m-1$, we have that for $t=2m-2$,
{\allowdisplaybreaks\begin{align*}
& \left \| M_{P_n}^{\circ}\delta_{2m-2} \right \|_{p}\\
&=\bigg(\sum_{1\leq i\leq m-1}(M_{P_n}^{\circ}\delta_{2m-2}(i))^p+\sum_{m\leq i\leq 2m-3}(M_{P_n}^{\circ}\delta_{2m-2}(i))^p\\
&\quad +(M_{P_n}^{\circ}\delta_{2m-2}(2m-2))^{p}+(M_{P_n}^{\circ}\delta_{2m-2}(2m-1))^{p}\bigg)^{1/p}\\
&=\left(m-1+(m-2)2^{-p}+1+1\right)^{1/p}\\
&=\left(m+1+2^{-p}(m-2)\right)^{1/p}.	
\end{align*}}When $1\leq k\leq m-2$, we have
{\allowdisplaybreaks\begin{align*}
\left \| M_{P_n}^{\circ}\delta_t \right \|_{p}&=\bigg(\sum_{1\leq i\leq k}(M_{P_n}^{\circ}\delta_{t}(i))^p+\sum_{k+1\leq i\leq 2k-1}(M_{P_n}^{\circ}\delta_{t}(i))^p+(M_{P_n}^{\circ}\delta_{t}(t))^{p}\\
&\quad +\sum_{2k+1\leq i\leq m+k-1}(M_{P_n}^{\circ}\delta_{t}(i))^p+\sum_{m+k\leq i\leq 2m-1}(M_{P_n}^{\circ}\delta_{t}(i))^p\bigg)^{1/p}\\
&=\left(k+(k-1)2^{-p}+1+(m-k-1)2^{-p}+(m-k)\right)^{1/p}\\
&=\left(m+1+2^{-p}(m-2)\right)^{1/p}.		
\end{align*}}

{\bf Case 2.} $n = 2m$ with $m \in \mathbb{N}_{+}$ and $m \geq 2$.

Following the approach from Case 1, we conclude that for any vertex $t \in V_{P_n}$,
\[
\left \| M_{P_n}^{\circ} \delta_t \right \|_{p} = \left( m + 1 + 2^{-p}(m-1) \right)^{1/p}.
\]

For $0<p\leq 1$, by Lemma $\ref{lem2-2}$, it is clear that for $m\geq 2$, we have
\begin{align*}
\left \| M_{P_n}^{\circ} \right \|_{\ell^p\to \ell^p}=\max_{t\in V_{P_n}}\left \| M_{P_n}^{\circ}\delta_t \right \|_{p}=
\begin{cases}
	\left(m+1+2^{-p}(m-2)\right)^{1/p} & n=2m-1,  \\
	\left(m+1+2^{-p}(m-1)\right)^{1/p}&  n=2m.
\end{cases}							
\end{align*}
						
{\rm(ii)} We now proceed to estimate the upper bound for $\left\Vert M_{P_n}^{\circ} \delta_t \right \Vert_{p}$ in the range $1 < p < \infty$.
Consider a function $f: V_{P_n} \to \mathbb{R}$ satisfying $\left\Vert f \right\Vert_p = 1$,
and let $t \in V_{P_n}$ be a vertex such that
\[
|f(t)|=\max_{v\in V_{P_n}}|f(v)|=:\eta.
\]
This maximality condition implies that for any vertex $i \in V_{P_n} \setminus \{t\}$,
we have the inequality
\begin{align}\label{equ4-7}
	|f(i)|^p\leq \min\{\eta^p,1-\eta^p\}.
\end{align}

Note that $M_{P_n}^{\circ}f(t)=|f(t)|$. For $ i \in V_{P_n}\setminus\{t\}$,
\begin{align*}
	\nonumber A_{P_{n},r}^{\circ}(|f|)(i) &=\frac{1}{|S(i,r)|}\sum_{y\in S(i,r)} |f(y)| \\
	&=
	\begin{cases}
		|f(i)| & \text{if } r=0,\\[5pt]
		\dfrac{|\tilde{f}(i-r)| + |\tilde{f}(i+r)|}{2} & \text{if } 1 \leq r \leq \min\lbrace i-1,\, n-i\rbrace,\\[5pt]
		|\tilde{f}(i-r)| + |\tilde{f}(i+r)| & \text{if } \min\lbrace i-1,\, n-i\rbrace < r \leq n-1.
	\end{cases}
\end{align*}
By Lemma \ref{lem4-2}, we have
\begin{align*}
	M_{P_n}^{\circ}f(i)
	&=\max\left\{\vert f(i)\vert,\max_{1\leq  r\leq \min\{i-1,n-i\}}\frac{|\tilde{f}(i-r)|+|\tilde{f}(i+r)|}{2},\right.\\
	&\quad \quad \quad \quad \quad \left.\max_{\min\{i-1,n-i\}<r\leq n-1}\left(|\tilde{f}(i-r)|+|\tilde{f}(i+r)|\right)\right\}.
\end{align*}

Clearly, the inequality
\begin{align}\label{equ4-8}
	M_{P_n}^{\circ}f(i)\leq \eta
\end{align}holds for any $i \in V_{P_n}$. Note that if $r=d_{P_n}(t,i)$ satisfies $\min\{i-1,n-i\}<r\leq n-1$, then
\begin{align*}
	|\tilde{f}(i-r)|+|\tilde{f}(i+r)|=\vert f(t)\vert,
\end{align*}
which implies that
\begin{align*}
	M_{P_n}^{\circ}f(i)=|f(t)|.
\end{align*}
Denote by $E$ the set of vertices $i \in V_{P_n}$ satisfying
\begin{align*} \min\{i-1,n-i\}<d_{P_n}(t,i)\leq n-1 \ \text{and} \  M_{P_n}^{\circ}f(i)=|f(t)|.
\end{align*}
We seek to understand how the set $E$ depends on the choice of the maximum vertex $t$, and to determine all vertices $i$ in $E$ for each $t$. The analysis requires separate cases depending on whether $n$ is odd or even.

{\bf Case 1.} $n = 2m-1$, where $m \in \mathbb{N}_{+}$ and $m \geq 2$.

Assume that $t$ is odd, i.e., $t = 2k - 1$ for some $k \in \mathbb{N}_+$ with $1 \leq k \leq m$.

If $2 \leq k \leq m - 1$, then it follows from
\[
\min\{i-1,n-i\}<d_{P_n}(t,i)\leq n-1
\]
that $1 \leq i \leq k - 1$ or $m + k \leq i \leq 2m - 1$. Thus for $2 \leq k \leq m - 1$, we have
\[
M_{P_n}^{\circ}f(i)=|f(t)|
\]
whenever $i = 1, 2, \ldots, k - 1, m + k, \ldots, 2m - 1$.

If $k = 1$, then $t = 1$. By the condition
\[
\min\{i-1,n-i\}\leq d_{P_n}(1,i)\leq n-1,
\]
we obtain $m + 1 \leq i \leq 2m - 1$. Therefore, for $k = 1$,
\[
M_{P_n}^{\circ}f(i)=|f(1)|
\]
for all $i = m + 1, m + 2, \ldots, 2m - 1$.

If $k = m$, then $t = n$. From
\[
\min\{i-1,n-i\}<d_{P_n}(n,i)\leq n-1
\]
it follows that $1 \leq i \leq m - 1$. Hence, for $k = m$,
\[
M_{P_n}^{\circ}f(i)=|f(2m-1)|
\]
for all $i = 1, 2, \ldots, m - 1$.

Now assume that $t$ is even, i.e., $t = 2k$ for some $k \in \mathbb{N}_+$ with $1 \leq k \leq m - 1$. Under the condition
\[
\min\{i-1,n-i\}<d_{P_n}(t,i)\leq n-1
\]
we conclude that
\begin{align*}
	M_{P_n}^{\circ}f(i)=|f(t)|
\end{align*}
for all $i = 1, 2, \ldots, k, m + k, \ldots, 2m - 1$.

Let $F:=E\bigcup \{t\}$. Clearly, for any vertex $i\in F$,  one has
\begin{align*}
	M^{\circ}_{P_n}f(i)=|f(t)|.
\end{align*}

Inequality \eqref{equ4-7} leads us to consider the following two cases: $\eta\leq 2^{-1/p}$ and $\eta>2^{-1/p}$.

If  $\eta\leq 2^{-1/p}$, then by \eqref{equ4-8}, we have
\begin{align*}
	\left\| M_{P_n}^{\circ}f\right\|_p =\left(\sum_{i=1}^{n}\vert M_{P_n}^{\circ}f(i)\vert^p\right)^{1/p}\leq (n/2)^{1/p}.
\end{align*}

If $\eta > 2^{-1/p}$, then the vertex $t$ is the unique maximizer. We now divide $V_{P_n}$ into $F$ and $V_{P_n}\setminus F$. For the vertex $i\in V_{P_n}\setminus F$, we will estimate the average function $A^{\circ}_{P_n,r}(|f|)$ according to different $r$.

For $r=0$, by \eqref{equ4-7}, we have
\[
A^{\circ}_{P_n,r}(|f|)(i)=|f(i)|\leq (1-\eta^p)^{1/p} \leq \frac{\eta+(1-\eta^p)^{1/p}}{2}.
\]

For $0 < r \leq \min\{i-1, n-i\}$, only one of the values $|\tilde{f}(i-r)|$ or $|\tilde{f}(i+r)|$ can be equal to $\eta$. Hence, by \eqref{equ4-7}, it follows that
\begin{equation*}
	A^{\circ}_{P_n,r}(|f|)(i) \leq \frac{\eta + (1 - \eta^p)^{1/p}}{2}.
\end{equation*}

For $\min\{i-1, n-i\} < r \leq n-1$, since $i \in V_{P_n} \setminus F$, the maximizer $t$ is not accessible. Hence, by \eqref{equ4-7}, we obtain
\[
A^{\circ}_{P_n,r}(|f|)(i) = |\tilde{f}(i-r)| + |\tilde{f}(i+r)| \leq (1-\eta^p)^{1/p} \leq \frac{\eta + (1-\eta^p)^{1/p}}{2}.
\]

Combining all the estimates above, we conclude that
\[
M^{\circ}_{P_n}f(x) = \max_{r=0,1,\ldots,n-1} A_{P_{n},r}^{\circ}(|f|)(x)\leq \frac{\eta+(1-\eta^p)^{1/p}}{2},
\]
and therefore
\begin{align*}
	\left\| M_{P_n}^{\circ}f\right\|_p&=\sum_{i\in F} \vert M_{P_n}^{\circ}f(i)\vert^p +\sum_{i\in V_{P_n}\setminus F} \vert M_{P_n}^{\circ}f(i)\vert^p\\
	\leq&\begin{cases}
		\left(m\eta^p+(m-1)(2^{-1}(\eta+(1-\eta^p)^{1/p}))^p\right)^{1/p}&  t=2k-1, \\
		\left((m+1)\eta^p+(m-2)(2^{-1}(\eta+(1-\eta^p)^{1/p}))^p\right)^{1/p}& t=2k.
	\end{cases}
\end{align*}

Noting that  $\eta > 2^{-1/p}$, we have for any $t \in V_{P_n}$ that
\begin{align*}
	\left\| M_{P_n}^{\circ}\right\|_{\ell^p\to \ell^p}\leq \left((m+1)\eta^p+(m-2)\left(\frac{1}{2}\left(\eta+(1-\eta^p)^{1/p}\right)\right)^p\right)^{1/p}.
\end{align*}						

{\bf Case 2.} $n = 2m$ with $m \in \mathbb{N}_{+}$ and $m \geq 2$.

The argument for this case proceeds exactly as in Case 1.

If $\eta\leq 2^{-1/p}$, then
\[
\left\| M_{P_n}^{\circ}\right\|_{\ell^p\to \ell^p}  \leq (n/2)^{1/p}.
\]

If $\eta> 2^{-1/p}$, then
\begin{align*}
\left\| M_{P_n}^{\circ}\right\|_{\ell^p\to \ell^p}\leq \left((m+1)\eta^p+(m-1)\left(\frac{1}{2}\left(\eta+(1-\eta^p)^{1/p}\right)\right)^p\right)^{1/p}.
\end{align*}

We continue to estimate the above equations, denote $a=\eta^p$ with $a\in \left(1/2,1\right]$, and set
\begin{align*}
h(a)=(m+1)a+(m-2)\left(\frac{1}{2}(a^{1/p}+(1-a)^{1/p})\right)^p,\\
\tilde{h}(a)=(m+1)a+(m-1)\left(\frac{1}{2}(a^{1/p}+(1-a)^{1/p})\right)^p.
\end{align*}

By Lemma \ref{lem4-3}, we have that
\begin{align*}
	h(a)\leq \frac{5}{4}h(1),\quad  \tilde{h}(a)\leq \frac{5}{4}\tilde{h}(1),
\end{align*}

We conclude that
\[
\left\| M_{P_n}^{\circ} \right\|_{\ell^p\to \ell^p} \leq
\begin{cases}
	\max\left\{ \left( \dfrac{5}{4} \left( m + 1 + \dfrac{m - 2}{2^p} \right) \right)^{1/p}, \left(\dfrac{2m - 1}{2}\right)^{1/p} \right\} & n = 2m - 1, \\[0.4cm]
	\max\left\{ \left( \dfrac{5}{4} \left( m + 1 + \dfrac{m - 1}{2^p} \right) \right)^{1/p}, m^{1/p} \right\} & n = 2m.
\end{cases}
\]
Clearly $m+1+2^{-p}(m-2)\geq (2m-1)/2$ and $m+1+(m-1)2^{-p}\geq m$, which proves {\rm(ii)}.

{\rm(iii)} For $p=\infty$, for any $t\in V_{P_n}$, we take $f=\delta_t$, then $\left\|M_{P_n}^{\circ}\delta_t\right\|_{\infty}=1$. Clearly $\left\|\delta_t\right\|_{\infty}=1$, thus
\[
\left\|M_{P_n}^{\circ}\right\|_{\infty}\geq 1.
\]
On the other hand, the equation $\eqref{equ4-8}$ confirms that
\[
\left\|M_{P_n}^{\circ}\right\|_{\infty}\leq 1,
\]
which proves {\rm(iii)}.

{\rm(iv)}
Assume $n=3$. Let $V_{P_3} = \{1, 2, 3\}$ and $E_{P_3} = \{1 \sim 2, 2 \sim 3\}$, and consider a function $f \colon V_{P_3} \to \mathbb{R}$ such that $\|f\|_p = 1$. The maximal function values are given by:
\begin{align*}
	M^{\circ}_{P_3}f(1) &= \max\left\{|f(1)|, |f(2)|, |f(3)|\right\}, \\
	M^{\circ}_{P_3}f(2) &= \max\left\{|f(2)|, \frac{|f(1)| + |f(3)|}{2}\right\}, \\
	M^{\circ}_{P_3}f(3) &= \max\left\{|f(1)|, |f(2)|, |f(3)|\right\}.
\end{align*}
We consider only the following three cases and the remaining cases can be handled analogously.

{\bf Case 1.} $|f(1)|\geq |f(2)|\geq |f(3)|$.

If $|f(2)|\geq (|f(1)|+|f(3)|)/2$, then we have
\begin{align*}
	M^{\circ}_{P_3}f(1)=M^{\circ}_{P_3}f(3)=|f(1)|, \quad 	M^{\circ}_{P_3}f(2)=|f(2)|.
\end{align*}
Thus
\begin{align*}
\left\|M^{\circ}_{P_3}f\right\|_p=(2|f(1)|^p+|f(2)|^p)^{1/p}\leq (3|f(1)|^p)^{1/p}\leq 3^{1/p}.
\end{align*}

If $|f(2)|< (|f(1)|+|f(3)|)/2$, then we have
\begin{align*}
	M^{\circ}_{P_3}f(1)=M^{\circ}_{P_3}f(3)=|f(1)|, \quad 	M^{\circ}_{P_3}f(2)=\frac{|f(1)|+|f(3)|}{2}.
\end{align*}
Thus
\begin{align*}
	\left\|M^{\circ}_{P_3}f\right\|_p=\left(2|f(1)|^p+2^{-p}(|f(1)|+|f(3)|)^p\right)^{1/p}
	\leq 3^{1/p}.
\end{align*}

{\bf Case 2.} $|f(1)|\geq |f(3)|\geq |f(2)|$.

We have  \begin{align*}
M^{\circ}_{P_3}f(1)=M^{\circ}_{P_3}f(3)=|f(1)|, \quad 	M^{\circ}_{P_3}f(2)=\frac{|f(1)|+|f(3)|}{2},
\end{align*}
and
\begin{align*}
\left\|M^{\circ}_{P_3}f\right\|_p=\left(2|f(1)|^p+2^{-p}(|f(1)|+|f(3)|)^p\right)^{1/p}\leq 3^{1/p}.
\end{align*}
						
{\bf Case 3.} $|f(2)|\geq |f(1)|\geq |f(3)|$.

We have
\begin{align*}	
M^{\circ}_{P_3}f(1)=M^{\circ}_{P_3}f(2)=M^{\circ}_{P_3}f(3)=|f(2)|,
\end{align*}
and
\begin{align*}
\left\|M^{\circ}_{P_3}f\right\|_p=(3|f(2)|^p)^{1/p}	\leq 3^{1/p}.
\end{align*}
						
Combining all cases, we conclude that for all $0<p<\infty$,
\[
\left\|M^{\circ}_{P_3}\right\|_{\ell^p\to \ell^p}\leq 3^{1/p}.
\]
Since the lower bound $\left\|M^{\circ}_{P_3}\right\|_{\ell^p\to \ell^p}\geq 3^{1/p}$ was established in part {(i)}, the proof is complete.
\end{proof}

\end{document}